\documentclass{elsarticle}

\usepackage[british]{babel}
\usepackage{amsmath, amsfonts, amsthm, amssymb}
\usepackage{multirow}
\usepackage{tikz}
		\usetikzlibrary{arrows, decorations.pathreplacing}
\usepackage{enumerate}
\usepackage{noindentafter}
\usepackage{algorithm}
\usepackage{algpseudocode}
\usepackage[colorlinks = true,linkcolor = blue, citecolor = blue]{hyperref}

\DeclareMathOperator*{\argmin}{arg\, min}
\DeclareMathOperator{\spn}{span}
\newcommand{\euler}{e}
\newcommand{\im}{i}
\newcommand{\real}{Re}
\newcommand{\mbbR}{\mathbb{R}}
\newcommand{\mbbC}{\mathbb{C}}
\newcommand{\df}{\mathrm{d}}
\newcommand\numberthis{\addtocounter{equation}{1}\tag{\theequation}}				
\newcommand{\hypref}[2]{\hyperref[#2]{#1 \ref*{#2}}}
\newcommand{\Break}{\State\textbf{break}}																	  
\newcommand{\wt}[1]{\widetilde{#1}}
\newcommand{\mc}[1]{\mathcal{#1}}
\newcommand{\rot}[1]{\makebox[1em][s]{\rotatebox{90}{#1}}}

\begin{document}


\begin{frontmatter}
		\title{Regula falsi based automatic regularization method for PDE constrained optimization}
		
		\author{Nick Schenkels\fnref{ref1}}
		\author{Wim Vanroose\fnref{ref2}}
		\address{Departement of Mathematics and Computer Science,\\University of Antwerp,\\ Antwerp, Belgium}
		\fntext[ref1]{Corresponding author: nick.schenkels@uantwerpen.be}
		\fntext[ref2]{wim.vanroose@uantwerpen.be}
		
		\begin{abstract}
				Many inverse problems can be described by a PDE model with unknown parameters
				that need to be calibrated based on measurements related to its solution.
				This can be seen as a constrained minimization problem where one wishes
				to minimize the mismatch between the observed data and the model predictions,
				including an extra regularization term, and use the PDE as a constraint.
				Often, a suitable regularization parameter is determined by solving the
				problem for a whole range of parameters -- e.g. using the L-curve -- which
				is computationally very expensive. In this paper we derive two methods
				that simultaneously solve the inverse problem and determine a suitable value
				for the regularization parameter. The first one is a direct generalization
				of the Generalized Arnoldi Tikhonov method for linear inverse problems.
				The second method is a novel method based on similar ideas, but with a
				number of advantages for nonlinear problems.
		\end{abstract}
		
		\begin{keyword}
				PDE constrained optimization, regularization, Morozov's discrepancy
				principle, Newton-Krylov, inverse scattering.
		\end{keyword}
\end{frontmatter}


\section{Introduction}
The dynamics of many complex applications are described by a PDE model $F(u, k) = 0$,
with solution or state variables $u$ and parameters or control variables $k$.
Examples include all forms of wave scattering problems \cite{abdoulaev2005, bruckner2017},
various financial models \cite{kaebe2009, inthout2010}, etc. The forward problem,
i.e. solving the PDE for $u$ given the parameters $k$, is often well understood
and is in many cases solved fast and accurately by a numerical method. The inverse
problem, i.e. finding the parameters $k$ such that the solution $u$ matches a
set of observations $\wt{u}$ as best as possible, is, however, much more complicated
because these problems are typically ill-posed.

Let $H(k)\in\mbbC^{m\times n}$, $u\in\mbbC^n$, $f(k)\in\mbbC^m$ and $k\in\mbbC^l$
be such that
\begin{equation}\label{eq:discretePDE}
		H(k)u = f(k),
\end{equation}
is the discretized version of the PDE with the appropriate initial and boundary
conditions. If $\wt{u}\in\mbbC^p$ are the observations of the solution of the PDE
and assuming that $u$ is an implicit function of $k$, we consider the following
constrained optimization problem:
\begin{equation}\label{eq:constropt}
		\left\{\begin{aligned}
				\min_{k\in\mbbC^l}\mc{J}(k) &= \min_{k\in\mbbC^l}\underbrace{\left\|Lu - \wt{u}
						\right\|^2}_{\mc{D}(k):=} + \alpha\underbrace{\left\|k - k_0
						\right\|^2}_{\mc{R}(k):=}\\
				H(k)u &= f(k).
		\end{aligned}\right.
\end{equation}
Here, $\left\|\cdot\right\|$ denotes the standard Euclidian norm, $L\in\mbbR^{p
\times n}$ is a linear operator that maps the full solution $u$ of the PDE to
the observed output $\wt{u}$, $\mc{D}(k)$ is a discrepancy or residual term
measuring the mismatch between the model predictions $u$ and the observations
$\wt{u}$ and $\mc{R}(k)$ is the regularization term added in order to  place
certain constraints on the parameters, incorporate prior knowledge, suppress
numerical errors or guarantee that the problem is well-posed.

There are now two difficulties, the first of which is calculating the gradient
of $\mc{J}(k)$. This is necessary because many nonlinear optimization algorithms
use some form of gradient information, e.g. Newton's method, steepest descent, nonlinear
CG, etc \cite{nocedal2006}. However, approximating $\nabla \mc{J}(k)$ using finite
difference methods is inefficient when, for example, $k$ is very high dimensional
\cite{tortorelli1994}. In order to avoid this, the adjoint method can be used in
order to calculate the gradient at the cost of only one PDE solve \cite{tortorelli1994,
plessix2006, jadamba2017}. The second difficulty is choosing the regularization
parameter $\alpha\in\mbbR^+$. Since this parameter models the balance between
model fidelity ($\alpha\rightarrow 0$) and the regularity of $k$ ($\alpha\rightarrow+\infty$),
its value greatly influences the reconstruction. While many papers describe how
the adjoint method can be used to solve nonlinear inverse problems, often the
regularization parameter is chosen by trial-and-error or using the L-curve
\cite{calvetti1999, calvetti2004, hansen2010, vogel2002}. Recent examples of
this include \cite{bruckner2017, jadamba2017}. These approaches, however, require
the solution of the inverse problem for many different values of the regularization
parameter, which is inefficient, computationally expensive and may take a long
time for large scale problems.

In this paper we derive two methods that simultaneously solve the inverse
problem and determine a suitable value for the regularization parameter. The
first method we call ``generalized Newton-Tikhonov'' (GNT) and is a direct
generalization of the Generalized Arnoldi Tikhonov method (GAT) for linear inverse
problems \cite{gazzola2014_2, gazzola2014, gazzola2015}. However, as we will
demonstrate, GNT requires a number of redundant computations which were not
needed in the original GAT method. The second method we call ``regula falsi
generalized Newton-Tikhonov'' (RFGNT) and is a novel method based on
similar ideas, but with a number of advantages which make it more efficient
for nonlinear problems.

The outline of the paper is as follows. In \hypref{section}{sec:lip} we give
an overview of how the regularization parameter can be chosen and how this is
automatically done by the generalized Arnoldi-Tikhonv method. This method will then
be used as a basis for the GNT and RFGNT algorithms we derive in \hypref{section}
{sec:nlip}. We then apply our methods to an inverse scattering problem in \hypref{section}
{sec:numexp}, where we use the adjoint method for the gradient computations, and
compare our results with other known regularization approaches.


\section{Automatic regularization for linear problems}\label{sec:lip}


\subsection{Choosing the regularization parameter}
It is well known that when dealing with inverse problems and measured data some
form of regularization is necessary in order to find a good solution \cite{hansen2010,
vogel2002}. However, the regularization parameter $\alpha$ can greatly influence
the outcome since it determines the balance between fitting the model to the noisy
data and the regularization. If, on the one hand, $\alpha$ is too small, the
regularization will have little to no effect and the noise in the data will corrupt
the outcome of the algorithm. If, on the other hand, $\alpha$ is too large, this
will lead to a solution that no longer fits the data very well. It may also
have lost many small details and be what is referred to as ``oversmoothed''.

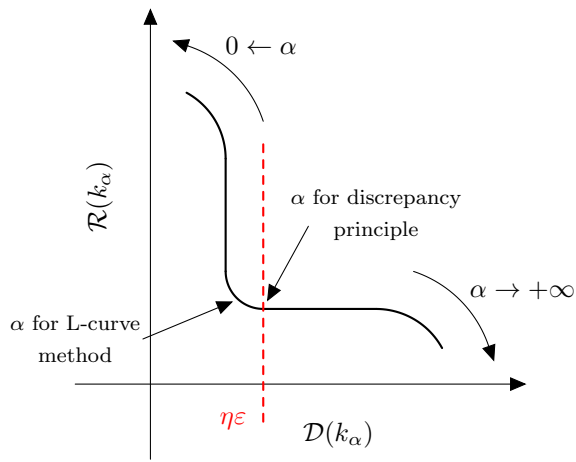
\begin{figure}[h]
		\centering
		\begin{tikzpicture}[line cap = round, line join = round, > = triangle 45]
		\draw[->] (-1, 0) -- (5, 0);
		\draw[->] (0, -1) -- (0, 5);
		
		 L-curve:
		\draw [shift = {(0, 3)}, thick]  plot[domain = 0:1.07,variable = \t] ({cos(\t r)},{sin(\t r)});
		\draw [shift = {(3, 0)}, thick]  plot[domain = 0.5:1.57,variable = \t] ({cos(\t r)},{sin(\t r)});
		\draw [shift = {(1.5, 1.5)}, thick]  plot[domain = 3.14:4.71,variable = \t] ({cos(\t r)/2},{sin(\t r)/2});
		\draw [thick] (1, 3) -- (1, 1.5);
		\draw [thick] (1.5, 1) -- (3, 1);
		
		\draw [dashed, thick, color = red] (1.5, -0.5) -- (1.5, 3.25);
		\node [anchor = north east] at (1.4, -0.25) {\color{red} $\eta\varepsilon$};
		
		\draw [shift = {(0, 3)}] plot[domain = 0.32:1.25, variable = \t] ({1*1.58*cos(\t r)+0*1.58*sin(\t r)},{0*1.58*cos(\t r)+1*1.58*sin(\t r)});
		\draw [shift = {(3, 0)}] plot[domain = 0.32:1.25, variable  =\t] ({1*1.58*cos(\t r)+0*1.58*sin(\t r)},{0*1.58*cos(\t r)+1*1.58*sin(\t r)});
		\draw [->] (0.5, 4.5) -- (0.29, 4.57);
		\draw [->] (4.5, 0.5) -- (4.57, 0.3);
		
		\node [label = below:$\mathcal{D}(k_\alpha)$] at (2.5, -0.25) {};
		\node [label = left:\rotatebox{90}{$\mathcal{R}(k_\alpha)$}] at (-0.25, 2.5) {};
		\node [font = \tiny, label = right:$0\leftarrow\alpha$] at (0.75, 4.5) {};
		\node [font = \tiny, label = right:$\alpha\rightarrow+\infty$] at (4, 1.25) {};
		
		
		\draw[->] (2.1, 2.1) -- (1.55, 1.1);
		\node [align = center] at (3, 2.25) {\footnotesize $\alpha$ for discrepancy\\ \footnotesize principle};
		\draw[->] (-0.1, 0.6) -- (1.1, 1.1);
		\node [align = center] at (-1, 0.6) {\footnotesize $\alpha$ for L-curve\\ \footnotesize method};
\end{tikzpicture}
		\caption{Sketch of the L-curve: the curve $(\mc{D}(k_\alpha), 
				\mc{R}(k_\alpha))$ typically has a rough L-shape. The L-curve
				method proposes to use the regularization parameter which
				corresponds to the corner of the L. The discrepancy principle on the
				other hand uses the value that corresponds to the intersection of the
				curve and the vertical line at $\eta\varepsilon$. This value is typically
				slightly bigger \cite{hansen1992}}
		\label{fig:LcurveDP}
\end{figure}

Let $k_\alpha$ be the solution of \eqref{eq:constropt} for a fixed regularization
parameter $\alpha$ and $u(k_\alpha)$ the corresponding solution to the PDE
\eqref{eq:discretePDE}. One way of determining a good value for $\alpha$ -- and
illustrating its effect on $k$ -- is the L-curve, see \hypref{figure}{fig:LcurveDP}.
By solving the inverse problem \eqref{eq:constropt} for a whole range of values
for $\alpha$ and looking at the the curve
\[
		(\mc{D}(k_\alpha), \mc{R}(k_\alpha)) = \left(\left\|Lu(k_\alpha) - \wt{u}
				\right\|^2, \left\|k_\alpha - k_0\right\|^2\right),
\]
it can be observed that it is roughly L-shaped. Heuristically, a ``good''
regularization parameter is the one that corresponds to the corner of the L,
since this will balance model fidelity and regularization \cite{hansen2010, vogel2002}.

Another way of choosing the regularization parameter $\alpha$ is the discrepancy
principle \cite{vogel2002, gazzola2014_2, morozov1984}, i.e. choose the regularization
parameter such that
\[
		\mc{D}(\alpha) := \mc{D}(k_\alpha) = \eta\underbrace{
				\left\|Lu - \wt{u}\right\|^2}_{\varepsilon:=}.
\]
Here, $\varepsilon$ is called the error norm and $1\leq\eta$ is a
tolerance value. The motivation behind this choice is that decreasing the
discrepancy $\mc{D}(\alpha)$ below the error norm will not necessarily
improve the reconstruction and can lead to overfitting. The downside of the 
discrepancy principle is that (an estimate of) $\varepsilon$ must be available.


\subsection{Generalized Arnoldi-Tikhonov}
The generalized Arnoldi-Tikhonov method was introduced in \cite{gazzola2014_2,
gazzola2014, gazzola2015} as a method to solve the classical Tikhonov problem for
linear problems of the form
\begin{equation}\label{eq:clasTik}
		\argmin_{x\in\mbbR^n}\left\|Ax - b\right\|^2 +
				\alpha\left\|x\right\|^2,
\end{equation}
with $x, b\in\mbbR^n$ and $A\in\mbbR^{n\times n}$. It is an iterative algorithm that
generates a sequence of approximations $x_0, x_1, x_2, \ldots$ that converge towards
the solution of \eqref{eq:clasTik}, while also updating the regularization parameter
in each iteration. This is done based on the discrepancy principle and the
current approximation of the solution. The method can best be understood by looking at the
discrepancy curve $(\alpha, \mc{D}(\alpha))$, see \hypref{figure}{fig:LcurveDPadj},
which can be seen as the analogue of the L-curve for the discrepancy principle.

\begin{figure}
		\centering
		\begin{tikzpicture}[line cap = round, line join = round, > = triangle 45]
		\draw[->] (-1, 0) -- (7, 0);
		\draw[->] (0, -1) -- (0, 5);
		
		\draw [thick] (0, 0.25) .. controls (2, 0.25) and (4, 2) .. (6, 4);
		
		\draw [dashed, thick, color = red] (-0.25, 2) -- (7, 2);
		\node [anchor = north] at (6, 2) {\color{red} $\eta\varepsilon$};
		
		\node [anchor = north] at (3.5, 0) {$\alpha$};
		\node [anchor = east] at (-0.1, 2.5) {\rotatebox{90}{$\mathcal{D}(\alpha)$}};
		
		\draw[->] (3.3, 2.6) -- (3.8, 2.1);
		\node [align = center] at (2.5, 2.75) {\footnotesize $\alpha$ for discrepancy\\  \footnotesize principle};
		
		\draw [thick, decorate, decoration={brace, mirror}] (0, -0.6) -- (3.7, -0.6)
				node [midway, yshift = -10] {\footnotesize Overfitting};
		\draw [thick, decorate, decoration={brace, mirror}] (3.8, -0.6) -- (7, -0.6)
				node [midway, yshift = -10] {\footnotesize Oversmoothing};
\end{tikzpicture}
		\caption{Plot of the discrepancy as a function of the regularization parameter.
				When $\alpha$ is small the model fidelity will be very high, but due to the
				noise in the data can exhibit overfitting. By increasing $\alpha$ more
				emphasis is put on the regularization term and overfitting is reduced.
				The discrepancy will start to increase however.}
		\label{fig:LcurveDPadj}
\end{figure}
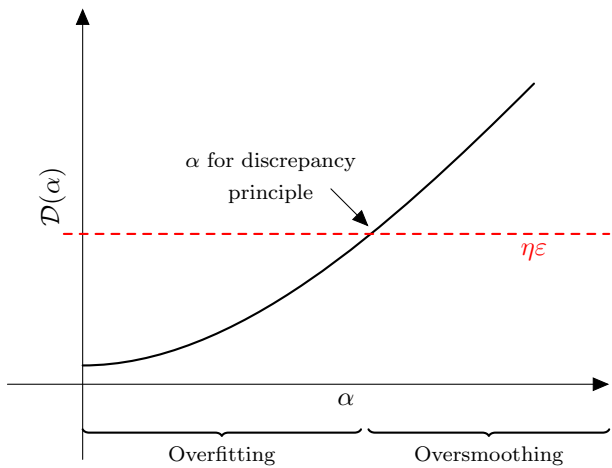

The idea behind GAT is to use the secant method in order to approximate the value of
$\alpha$ for which $\mc{D}(\alpha) = \eta\varepsilon$. The method is also a Krylov
subspace method based on the Arnoldi decomposition of the matrix $A$ \cite{saad2003,
vandervorst2003}. This means that the iterates for the solution of \eqref{eq:clasTik}
are given by
\[
		x_{\alpha, i} := \argmin_{x\in\mc{K}_i}\left\|Ax - b\right\|^2 + \alpha\left\|x\right\|^2,
\]
with
\[
		\mc{K}_i = \mc{K}_i(A, b) = \spn\left\{b, Ab, A^2b, \ldots, A^{i - 1}b\right\}
\]
the associated Krylov subspace of dimension $i$. In each iteration, a new basis vector
is added to the Krylov subspace and the iterates are updated in order to account for
this new basis vector. It is important to note that the constructed Krylov basis is independent
of the regularization parameter. This means that it can be stored and reused in the next
iteration when the regularization parameter is updated. In order to update the
current best estimate for the regularization parameter $\alpha_{i - 1}$, GAT assumes that
\eqref{eq:clasTik} is simultaneously solved without regularization, i.e. for $\alpha = 0$.
This means that the points $(0, \mc{D}(0))$ and $(\alpha_{i - 1}, \mc{D}(\alpha_{i - 1}))$
on the discrepancy curve are known and
the regularization parameter can be updated using one step of the secant method:
\[
		\alpha_i = \frac{\eta\varepsilon - \mc{D}(0)}{\mc{D}(\alpha_{i - 1}) - \mc{D}(0)}\alpha_{i - 1}
\]
Furthermore, instead of solving \eqref{eq:clasTik} to convergence each time the regularization
parameter is updated and calculate the value $\mc{D}(\alpha_{i - 1})$ exactly,
the inverse problem is solved in the currently constructed Krylov subspace. This means
that if
\[
		\mc{D}_i(\alpha) := \mc{D}(x_{\alpha, i}) = \left\|Ax_{\alpha, i} - b\right\|^2
\]
is the discrepancy after $i$ iterations -- or equivalently the discrepancy in the Krylov subspace $\mc{K}_i$
of dimension $i$ -- then the GAT update for the regularization
parameter is given by:
\begin{equation}\label{eq:alphaupdate}
		\alpha_i = \left|\frac{\eta\varepsilon - \mc{D}_i(0)}{\mc{D}_i(\alpha_{i - 1}) - \mc{D}_i(0)}\right|\alpha_{i - 1}.
\end{equation}
Then, once the discrepancy principle is satisfied, i.e. $\mc{D}_i(\alpha_{i - 1}) \leq\eta\varepsilon$,
the algorithm is stopped.

Note that the absolute value is added because now it is possible for both $\mc{D}_i(0)$
and $\mc{D}_i(\alpha_{i - 1})$ to be bigger than $\eta\varepsilon$, which can otherwise
result in negative values for the regularization parameter. This only happens in
the first few iterations when the constructed Krylov subspace is to small to contain
a good approximation for the solution. Since the constructed Krylov basis is
independent of $\alpha$, the fact that the regularization parameter is estimated
incorrectly in the first few iterations does not matter. It is typically
only when $\mc{D}_i(0)$ becomes smaller than $\eta\varepsilon$ that the estimates
start to improve. Finally, we remark that the value $\alpha = 0$ is chosen and fixed
for the secant updates because in this case the linear system that needs to be solved
is smaller and the method becomes equivalent to the GMRES algorithm \cite{saad1986}. 
However, since the algorithm is stopped once the discrepancy principle is satisfied,
this may result in an underestimation of the regularization parameter.


\section{Automatic regularization for nonlinear problems}\label{sec:nlip}


\subsection{Generalized Newton-Tikhonov}
Although GAT is a Krylov subspace method for linear the linear Tikhonov problem, 
the same idea can be used for nonlinear problems. The update for the regularization
parameter \eqref{eq:alphaupdate} can even be applied directly if we change the notation
back to our nonlinear problem. The only difference is how the iterates are calculated.
If we use Newton's method to solve \eqref{eq:constropt}
and $k_{\alpha, i}$ is the $i$th Newton iteration using a fixed regularization parameter
and $u(k_{\alpha, i})$ the corresponding solution to the PDE \eqref{eq:discretePDE}, then
the discrepancy after $i$ iterations is now given by
\begin{equation}\label{eq:nldiscrp}
		\mc{D}_i(\alpha) := \mc{D}(k_{\alpha, i}) = \left\|Lu(k_{\alpha, i}) - \wt{u}\right\|^2.
\end{equation}
An overview of this method, which we will call generalized Newton-Tikhonov (GNT),
is given in \hypref{algorithm}{alg:gnt}. This method should be seen as a direct
generalization of the GAT algorithm for nonlinear problems. As we will demonstrate
with our numerical experiments, the method can be used to solve our nonlinear inverse problem,
but it has a number of drawbacks that were not present in the original GAT method
for linear inverse problems.

\begin{algorithm}
		\caption{generalized Newton-Tikhonov (GNT)}\label{alg:gnt}
		\begin{algorithmic}[1]
				\State Choose initial $\alpha_0$, $k_0$
				\For{$i = 1, \ldots,$ maxIter}
						\State Calculate $k_{0, i}$ based on $k_{0, i - 1}$.
								\Comment{This requires 1 Newton step.$\:$}\label{alg:gnt:3}
						\State Calculate $k_{\alpha_{i - 1}, i}$ based on $k_0$.
								\Comment{This requires i Newton steps.}\label{alg:gnt:4}
						\State Calculate $\mc{D}_i(0)$ and $\mc{D}_i(\alpha_{i - 1})$ using \eqref{eq:nldiscrp}.
						\If{$\mc{D}_i(\alpha_{i - 1})\leq\eta\varepsilon$}
								\Break
						\Else
								\State Calculate $\alpha_i$ using \eqref{eq:alphaupdate}.\label{alg:gnt:9}
						\EndIf
				\EndFor
		\end{algorithmic}
\end{algorithm}

A first thing to note is that both GAT and GNT update the regularization parameter
based on the regularized and the non-regularized solution after a certain number of
iterations. For the non-regularized iterations, this means that in each GNT iteration
we need to perform one Newton step, see \hypref{algorithm}{alg:gnt} \hypref{line}{alg:gnt:3}.
However, when we wish to determine $k_{\alpha_{i - 1}, i}$, we cannot use previous
best approximation $k_{\alpha_{i - 2}, i - 1}$. This is
because the Newton iterations depend on $\alpha$, so we have to restart them from
$k_0$ and perform $i$ new Newton steps, see \hypref{algorithm}{alg:gnt} \hypref{line}{alg:gnt:4}. In the original GAT method for linear
problems this is not an issue, since the Krylov basis is independent from
the regularization parameter. Therefore, in each iteration only one new Krylov
basis vector has to be determined and added to the current basis. For GNT on the
other hand, in each iteration the regularized Newton iterations have to be restarted
and the number of Newton steps that needs to be computed increases each time. In order to 
to perform $i$ GNT iterations, the number of Newton iterations needed is:
\[
		\underbrace{(1 + 1 + \ldots + 1)}_{i\text{ times \hypref{line}{alg:gnt:3}}} +
				\underbrace{(1 + 2 + \ldots + i)}_{i\text{ times \hypref{line}{alg:gnt:4}}}
				= \frac{i(i + 3)}{2}.
\]

Another thing to remark is that in the original GAT method the reason
for using $\alpha = 0$ as the second point for the secant step and not updating it,
is because for this choice of $\alpha$ the linear system that needs to be solved
is smaller and hence, easier to solve. To draw the parallel with the original method
we also use this value for GNT. However, for nonlinear problems there is no direct
benefit for using this value and if a better initial estimate for the regularization
parameter is available or if $\alpha$ must be larger than $0$ in order for the
inverse problem to be well posed, another value can be used. This might also improve
the quality of the estimate for the regularization parameter.


\subsection{Regula falsi generalized Newton-Tikhonov}
Two drawbacks of GNT are the increasing number of Newton iterations and the 
fact that the value for the regularization parameter only slowly converges to
the value satisfying the discrepancy principle, which we will refer to as $\alpha^*$.
This can be seen in our numerical experiments, see \hypref{figure}{fig:details}, and
is explained by the fact that the secant update step for the regularization parameter
is done with a fixed point at $\alpha = 0$ and the discrepancy is only calculated
up to a limited number of Newton iterations, see \hypref{algorithm}{alg:gnt}
\hypref{line}{alg:gnt:9}. Therefore, in order to limit the total number of Newton iterations
and better approximate $\alpha^*$, we propose an alternative approach based on the
regula falsi method.

Assuming we have values $\alpha_0, \alpha_1\in\mbbR^+$ such that $\alpha_0\leq\alpha^*\leq\alpha_1$,
we determine the line between $(\alpha_0, \mc{D}(\alpha_0)$ and $(\alpha_1, \mc{D}(\alpha_1))$
and take $\alpha_2$ as the value for which this equals the discrepancy:
\begin{equation}\label{eq:alphaupdate2}
		\alpha_2 = \frac{\eta\varepsilon - \mc{D}(\alpha_0)}{
				\mc{D}(\alpha_1) -	\mc{D}(\alpha_0)}\left(\alpha_1 -
				\alpha_0\right) + \alpha_0
\end{equation}
We then solve \eqref{eq:constropt} with $\alpha_2$ to determine
$\mc{D}(\alpha_2)$ and replace either $\alpha_0$ or $\alpha_1$ with $\alpha_2$,
such that the interval $[\alpha_0, \alpha_1]$ will always contain $\alpha^*$.
Furthermore, in contrast to GNT, we will calculate $\mc{D}(\alpha)$ and
$k_\alpha$ exactly and not up to a limited number of iterations. This is justified
in the GAT algorithm, where due to the presence of a basis of the Krylov subspace
the discrepancy after $i$ iterations could be acquired at a low cost. However, as
we saw with GNT, this is no longer the case for nonlinear problems. In order to
limit the total number of required Newton iterations we will instead update the initial
guess for the Newton iterations every time $\alpha_2$ is updated. Similarly to the
regula falsi step for the regularization parameter, we take a weighted
linear approximation based on the solutions in $\alpha_0$ and $\alpha_1$, i.e.
$k_{\alpha_0}$ and $k_{\alpha_1}$, which were already calculated:
\begin{equation}\label{eq:k0update}
		k_{\alpha_2, 0} = \frac{\alpha_2 - \alpha_0}{\alpha_1 - \alpha_0}\left(
				k_{\alpha_1} - k_{\alpha_0}\right) + k_{\alpha_0}
\end{equation}
Then, once the value for $\alpha_2$ starts to converge, we terminate the algorithm,
see \hypref{algorithm}{alg:rfgnt} \hypref{line}{alg:rfgnt:9}.

An overview of this algorithm, which we call regula falsi generalized Newton-Tikhonov
(RFGNT), is given in \hypref{algorithm}{alg:rfgnt}. It should be noted that the only
difference between the updates for the regularization parameter in GNT and in RFGNT
is that in RFGNT $\alpha_0$ does not have the fixed value $0$. This also means that
when the initial interval $[\alpha_0, \alpha_1]$ does not contain $\alpha^*$, we
can update the interval based on the secant method until it does, i.e. $\alpha_0\leftarrow\alpha_1$
and $\alpha_1\leftarrow\alpha_2$ (or vise versa if $\alpha_2 < \alpha_1$).

\begin{algorithm}
		\caption{Regula falsi generalized Newton-Tikhonov (RFGNT)}\label{alg:rfgnt}
		\begin{algorithmic}[1]
				\State Choose initial $k_0$.
				\State Choose initial $\alpha_0$ and $\alpha_1$ such that $\alpha_0\leq
						\alpha^*\leq\alpha_1$.
				\State Calculate $k_{\alpha_0}$ and $\mc{D}(\alpha_0)$ by
						solving \eqref{eq:constropt} (starting with initial $k_0$).\label{alg:rfgnt:1}
				\State Calculate $k_{\alpha_1}$ and $\mc{D}(\alpha_1)$ by
						solving \eqref{eq:constropt} (starting with initial $k_0$).\label{alg:rfgnt:2}
				\For{$i = 1, \ldots,$ maxIter}
						\State Calculate $\alpha_2$ using \eqref{eq:alphaupdate2}.
						\State Calculate $k_{\alpha_2, 0}$ using \eqref{eq:k0update}.
						\State Calculate $k_{\alpha_2}$ and $\mc{D}(\alpha_2)$ by
								solving \eqref{eq:constropt} (starting with initial $k_{\alpha_2, 0}$).\label{alg:rfgnt:8}
						\If{$\left|\alpha_2 - \alpha_2^{old}\right|/\alpha_2^{old} < 10^{-3}$}\label{alg:rfgnt:9}
								\Break
						\Else
								\State $\alpha_2^{old}\leftarrow\alpha_2$
								\State Replace $\alpha_0$ or $\alpha_1$ with $\alpha_2$
										based on $\mc{D}(\alpha_2)$.
						\EndIf
				\EndFor
		\end{algorithmic}
\end{algorithm}


\section{Numerical experiments}\label{sec:numexp}


\subsection{Inverse scattering}
Consider the homogeneous Helmholtz equation
\[
		\left(\Delta + k^2\right)u_{tot} = 0.
\]
on a square domain $\Omega\subseteq\mbbR^2$ with exterior complex scaling (ECS) boundary
conditions and a spatially varying wave number $k:\Omega\longrightarrow\mbbR$.
We will assume that the total wave $u_{tot}:\Omega\longrightarrow\mbbC$ can be written
as the sum of an incoming wave and the resulting scattered wave:
\[
		u_{tot} = u_{in} + u_{sc}.
\]
If for multiple incoming waves of the form
\[
		u_{in}^\theta(x, y) = \euler^{ik_0(\cos\theta x + \sin\theta y)},
\]
the resulting scattered waves are given by $u_{sc}^\theta$, then this can be written
as one big system of equations:
\[
		\underbrace{\begin{pmatrix}\left(\Delta + k^2\right)\\&\left(\Delta + k^2\right)
				\\&&\ddots\\&&&\left(\Delta + k^2\right)\end{pmatrix}}_{\mc{H}:=}\begin{pmatrix}
				u_{in}^{\theta_1} + u_{sc}^{\theta_1}\\u_{in}^{\theta_2} + u_{sc}^{\theta_2}\\
				\vdots\\u_{in}^{\theta_t} + u_{sc}^{\theta_t}\end{pmatrix} = 
				\begin{pmatrix}0\\0\\\vdots\\0\end{pmatrix}
\]
If we denote $u_{in} = \left(u_{in}^{\theta_1}, \ldots, u_{in}^{\theta_t}\right)^T$
and $u_{sc} = \left(u_{sc}^{\theta_1}, \ldots, u_{sc}^{\theta_t}\right)^T$, then
this is equivalent to
\begin{equation}\label{eq:constr}
		\mc{H}u_{sc} = \underbrace{\left(k_0^2 - k^2\right)u_{in}}_{\mc{F}:=}.
\end{equation}
For our numerical experiment, we will try to reconstruct the wave number $k$ based
on measurements of the scattered wave at the boundary $\delta\Omega$. If $Hu = f$
is the discrete version of \eqref{eq:constr}, $\wt{u}$ the measured values at $\delta\Omega$ and $L$
the restriction operator that maps the full solution $u$ on the discretized domain
to its values on boundary of the domain, then we get the following constrained minimization
problem:
\begin{equation}
		\left\{\begin{aligned}
				\min_{k\in\mbbR^n}\mc{J}(k) &= \min_{k\in\mbbR^n}\left\|Lu - \wt{u}
						\right\|_2^2 + \alpha\left\|k - k_0\right\|_2^2\\
				Hu &= f
		\end{aligned}\right.
\end{equation}


\subsection{The adjoint method}\label{sec:adj}
Because we will use Newton's method for the optimization in GNT and RFGNT, we
calculate the gradient of $\mc{J}(k)$ using the adjoint method. Let $\left<x\mid y\right>
= x^*y$ denote the complex inner product for $x, y\in\mbbC^n$ and $x^*$ the conjugate
transpose of $x$, then we can write the cost function as
\[
		\mc{J}(k) = \left<Lu - \wt{u}\left|Lu - \wt{u}\right.\right>
				+ \alpha\left<k - k_0\left|k - k_0\right.\right>.
\]
It follows that:
\begin{align*}
		\frac{\df\mc{J}}{\df k}\ =&\ \left<Lu - \wt{u}\left| L\frac{\df u}{\df k}
				\right.\right> + \left<\left.L\frac{\df u}{\df k}\right|Lu - \wt{u}\right>\\
		& +\alpha\left<1\left|k - k_0\right.\right> + \alpha\left<\left.k - k_0\right|1\right>\\
		=&\ 2\real\left(\left<L^*(Lu - \wt{u})\left|\frac{\df u}{\df k}\right.
				\right> + \alpha\left<\left.k - k_0\right|1\right>\right).
\end{align*}
The difficulty is now the derivative of the state variables with respect to the control
variables, i.e. $\df u/\df k$. In order to avoid needing to calculate this term, the
adjoint method introduces an adjoint variable $\lambda\in\mbbC^m$ as the solution of
\begin{equation}\label{eq:adjeq}
		H^*(k)\lambda = L^*\left(Lu - \wt{u}\right),
\end{equation}
where $H$ was the matrix representing the discretized PDE, see \eqref{eq:discretePDE}.
From \eqref{eq:discretePDE} it also follows that
\[
		\frac{\df H}{\df k}u + H\frac{\df u}{\df k} =	\frac{\df f}{\df k}.
\]
Substituting this into the derivative of $\mc{J}$, we can eliminate the term
$\df u/\df k$:
\begin{align*}
		\frac{\df\mc{J}}{\df k} &= 2\real\left(\left<H^*\lambda\left|\frac{\df u}
				{\df k}\right.\right> + \alpha\left<\left.k - k_0\right|1\right>\right)\\
		&= 2\real\left(\left<\lambda\left|H\frac{\df u}{\df k}\right.\right>
				+ \alpha\left<\left.k - k_0\right|1\right>\right)\\
		&= 2\real\left(\left<\lambda\left|\frac{\df f}{\df k} - \frac{\df H}{\df k}u
				\right.\right> + \alpha\left<\left.k + k_0\right|1\right>\right).
				\numberthis\label{eq:djdk}
\end{align*}
This means that the gradient of the cost function can now be evaluated in three steps:
\begin{enumerate}[i)]
		\item Given $k$, solve the original PDE, see \eqref{eq:discretePDE}, for $u$.
		\item Use $u$ to solve the adjoint PDE, see \eqref{eq:adjeq}, for $\lambda$.
		\item Evaluate the gradient $\df\mc{J}/\df k$ using \eqref{eq:djdk}.
\end{enumerate}
Note that solving the adjoint equation \eqref{eq:adjeq} is closely related to the
original PDE \eqref{eq:discretePDE} and that the cost of solving it will be similar.
This also means that we can evaluate the gradient of $\mc{J}(k)$ at the cost of
only one extra PDE solve per iteration. For more information on the adjoint method
we refer to \cite{tortorelli1994, plessix2006}.


\subsection{Newton-Krylov}
In order to solve the inverse problem \eqref{eq:constropt}, we will use an algorithm
based on the line search Newton-CG method described in \cite{nocedal2006}. This is
a Newton-Krylov method, where the parameters $k$ are updated using Newton's method, i.e.
\[
		k_{\alpha, i + 1} = k_{\alpha, i} + \gamma\Delta k,
\]
and the Hessian system for the Newton search direction $\Delta k$ is solved using a Krylov
subspace method, CG in this case \cite{nocedal2006, shewchuk1994}:
\[
		\nabla^2\mc{J}(k_{\alpha, i})\Delta k= -\nabla\mc{J}(k_{\alpha, i}).
\]
The CG iterations are stopped once the residual is smaller than
\[
		\min\left(0.5, \sqrt{\left\|\nabla\mc{J}(k_{\alpha, i})\right\|}\right)\nabla\mc{J}(k_{\alpha, i})
\]
or $\left\|\Delta k\right\| < 10^{-3}$. The line search will be a simple backtracking algorithm
starting from $\gamma = 1$ and halving this value until
\[
		\mc{J}(k_{\alpha, i} + \gamma\Delta k) < \mc{J}(k_{\alpha, i}).
\]
We terminate the Newton iterations once they start to stagnate, i.e.
\[
		\left|\frac{J(k_{\alpha, i + 1}) - J(k_{\alpha, i})}{J(k_{\alpha, i})}\right| < 10^{-3}
\]

It should be noted that each function evaluation comes at the cost of one PDE solve. Furthermore,
each CG iteration requires one matrix vector product with the Hessian. In order to
avoid this, a finite difference approximation can be used. Using a central difference
scheme the approximation is given by:
\begin{equation}\label{eq:hessapprox}
		\nabla^2\mc{J}(k)v\approx\frac{\nabla\mc{J}(k + hv) - \nabla\mc{J}(k - hv)}{2h}.
\end{equation}
This implies that each Newton iterations requires three gradient calculations.
Using the adjoint method, this is equivalent to three solves of the original PDE and three
solves of the adjoint PDE.  By replacing the central difference approximation of the Hessian matrix-vector product
 \eqref{eq:hessapprox} with a forward or backward difference approximation, this can
be reduced to two. It should be noted that the Hessian matrix-vector product can also be determined using the
second order adjoint method, see \cite{tortorelli1994, jadamba2017, wang1992}. For
simplicity we chose not to do so.


\subsection{Discretization}
The discretization of the Helmholtz systems in $\mc{H}$ is independent with respect
to the angle of the incoming wave. $H$ will therefore be a block diagonal matrix
with each block being a discrete version of the operator $(\Delta + k^2)$. For our
numerical experiment we take $\Omega = [-5, 5]^2$ and discretize it using a regular
$200\times 200$ grid with grid spacing $h$. We add a small buffer zone of 10 grid points
before the points where the measurements are taken and then add another 10 grid
points before the start of the complex tails for the exterior complex scaling. 

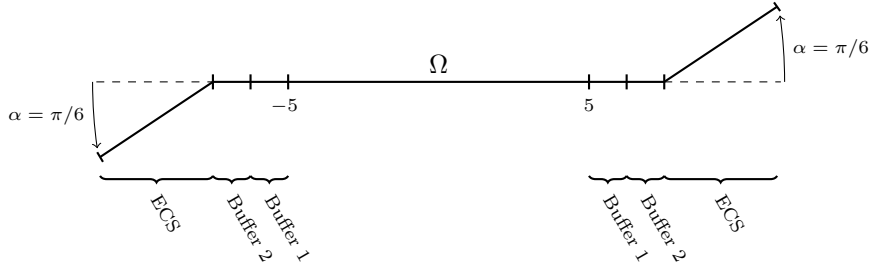
\begin{figure}
		\centering
		\begin{tikzpicture}[xscale = 0.5]
		\draw[thick] (-6, 0) -- (6, 0);
		\draw[thick] (-6, -0.1) -- (-6, 0.1);
		\draw[thick] (-5, -0.1) -- (-5, 0.1);
		\draw[thick] (-4, -0.1) -- (-4, 0.1);
		\draw[thick] (4, -0.1) -- (4, 0.1);
		\draw[thick] (5, -0.1) -- (5, 0.1);
		\draw[thick] (6, -0.1) -- (6, 0.1);
		
		\node (A) at (-9, -1) {};
		\draw[thick] (-9, -1) -- (-6, 0);
		\begin{scope}[rotate around = {50:(A)}]
				\draw[thick] (-9, -1.1) -- (-9, -0.9);
		\end{scope}
		\node (B) at (9, 1) {};
		\draw[thick] (6, 0) -- (9, 1);
		\begin{scope}[rotate around = {55:(B)}]
				\draw[thick] (9, 1.1) -- (9, 0.9);
		\end{scope}
		\draw[dashed] (-9.1, 0) -- (-6, 0);
		\draw[->] (-9.2, 0) arc (180:196:3.2) node[midway, left]{\scriptsize $\alpha = \pi/6$};
		\draw[dashed] (6, 0) -- (9.1, 0);
		\draw[->] (9.2, 0) arc (0:16:3.2)  node[midway, right]{\scriptsize $\alpha = \pi/6$};
		
		\node[above] at (0, 0) {$\Omega$};
		\node[below] at (-4.1, -0.1) {\scriptsize $-5$};
		\node[below] at (4, -0.1) {\scriptsize $5$};
		\draw[thick, decorate, decoration = {brace, mirror}] (-9, -1.25) -- (-6, -1.25);
		\draw[thick, decorate, decoration = {brace, mirror}] (-6, -1.25) -- (-5, -1.25);
		\draw[thick, decorate, decoration = {brace, mirror}] (-5, -1.25) -- (-4, -1.25);
		\draw[thick, decorate, decoration = {brace, mirror}] (4, -1.25) -- (5, -1.25);
		\draw[thick, decorate, decoration = {brace, mirror}] (5, -1.25) -- (6, -1.25);
		\draw[thick, decorate, decoration = {brace, mirror}] (6, -1.25) -- (9, -1.25);
		\node[below] at (-7.25, -1.35) {\rotatebox{-60}{\scriptsize ECS}};
		\node[below] at (-5, -1.35) {\rotatebox{-60}{\scriptsize Buffer 2}};
		\node[below] at (-4, -1.35) {\rotatebox{-60}{\scriptsize Buffer 1}};
		\node[below] at (5, -1.35) {\rotatebox{-60}{\scriptsize Buffer 1}};
		\node[below] at (6, -1.35) {\rotatebox{-60}{\scriptsize Buffer 2}};
		\node[below] at (7.75, -1.35) {\rotatebox{-60}{\scriptsize ECS}};
\end{tikzpicture}
		\caption{Exterior complex scaling in 1D: by adding complex tails to the domain
				and Dirichlet boundary conditions at the end of these tails, outgoing waves
				at the boundary $\partial\Omega$ can be simulated. The angle with respect
				to the real axis and how far the complex tails extend into the complex plain
				can be chosen.}
		\label{fig:ECS}
\end{figure}

The idea behind exterior complex scaling \cite{mccurdy2004} is to extend the domain
into the complex plane. By imposing Dirichlet boundary conditions at the end of
these complex ``tails'', outgoing waves at the boundary $\delta\Omega$ can be simulated.
Numerically we do this by adding points to the real domain (with the same spacing $h$)
and rotating them into the complex plane under a chosen angle, see \hypref{figure}{fig:ECS}.
We use 80 grid points for the complex tails -- one third of the real domain --
in each direction with an angle $\alpha = \pi/6$. This means that the full grid has size
$400\times 400$.

To discretize the Laplace operator we use second order finite differences, which on a
regular grid with spacing $h$, is given by the formula:
\begin{equation}\label{eq:discretelap}
		\begin{aligned}
				\Delta u(x, y)\approx&\phantom{+}\frac{u(x - h, y) - 2u(x, y) +	u(x + h, y)}{h^2}\\
																		&+\frac{u(x, y - h) - 2u(x, y) +  u(x, y + h)}{h^2}
		\end{aligned}
\end{equation}
However, because of the way the complex tails are constructed, the full grid is no
longer regular. Therefore, the general form for irregular grids has to be used.
Furthermore, we will assume that the wave number $k\in\mbbR^n$ with $n = 200^2$
is equal to a base value $k_0 = 1$ outside of $\Omega$. Inside the region of interest
itself we add an offset based on the sum of three Gaussian functions placed symmetrically
on a circle with radius $r = 2.5$:
\[
		k = k_0\sqrt{1 + \chi},
\]
with
\begin{align*}
		\chi(x, y) =&\ \euler^{-(x - r)^2 - y^2} + \euler^{-\left(x - r\cos\left(\frac{2\pi}{3}
				\right)\right)^2 - \left(y - r\sin\left(\frac{2\pi}{3}\right)\right)^2}\\
		&+\ \euler^{-\left(x - r\cos\left(\frac{4\pi}{3}\right)\right)^2 - \left(y -
				r\sin\left(\frac{4\pi}{3}\right)\right)^2}.
\end{align*}

In order to simulate the measurements, we generate incoming waves of the form
\[
		u_{in}^\theta(x, y) = \euler^{ik_0(\cos\theta x + \sin\theta y)}
\]
for $50$ different value of $\theta\in[0, 2\pi[$ and calculate the corresponding
scattered waves $u_{exact}\in\mbbC^{50\cdot 400^2}$. Taking into account the $10$ buffer
points we added, we have $876$ observations for every angle $\theta$, which we
select using the matrix $L\in\mbbR^{50\cdot 876\times 50\cdot 400^2}$, see \hypref{figure}{fig:scatwave}.
We then add random Gaussian white noise to generate the measurements
\[
		\wt{u} = Lu_{exact} + \sigma\left(e_1 + \im e_2\right)
\]
with $e_1$ and $e_2\sim\mc{N}(0, I_{50\cdot 876})$ and choose $\sigma$ such that the noise level
is approximately $10\%$:
\[
		\frac{\left\|Lu_{exact} - \wt{u}\right\|^2}
				{\left\|Lu_{exact}\right\|^2}\approx 0.10
\]
This resulted in a value $\sigma = 0.075$ and $\varepsilon = 492.7031$,
which we use for the discrepancy principle combined with $\eta = 1$. Because of the
diagonal structure of \eqref{eq:constr} we also do not construct the matrix $H$ explicitly,
but rather solve the $50$ Helmholtz systems for the different angles in parallel,
where each discrete version of the Helmholtz system $(\Delta + k^2)$ has size
$160000\times 160000$.

\begin{figure}
		\centering
		\includegraphics[width = 0.33\linewidth]{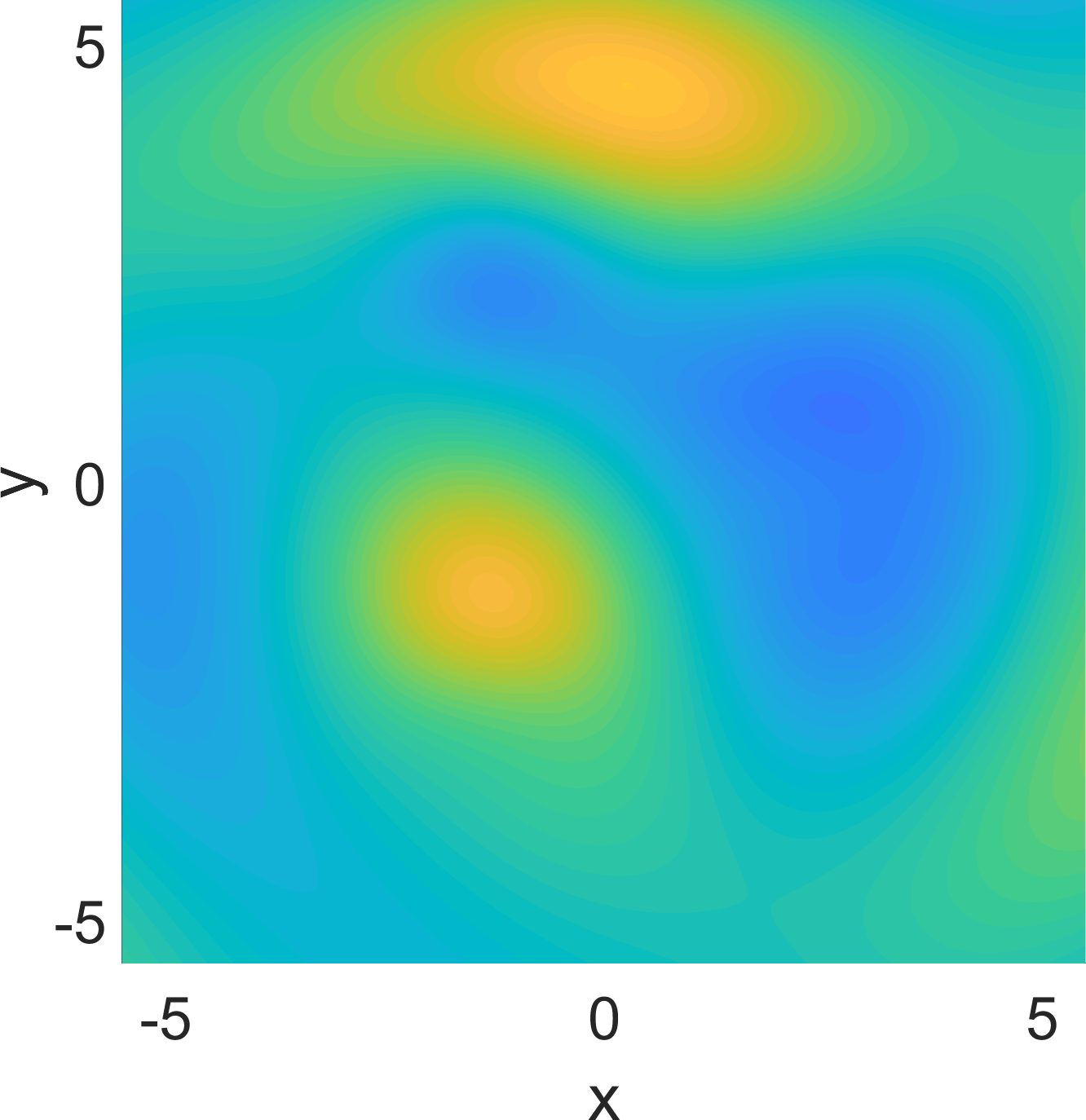}\hspace{2.5pt}
		\includegraphics[width = 0.33\linewidth]{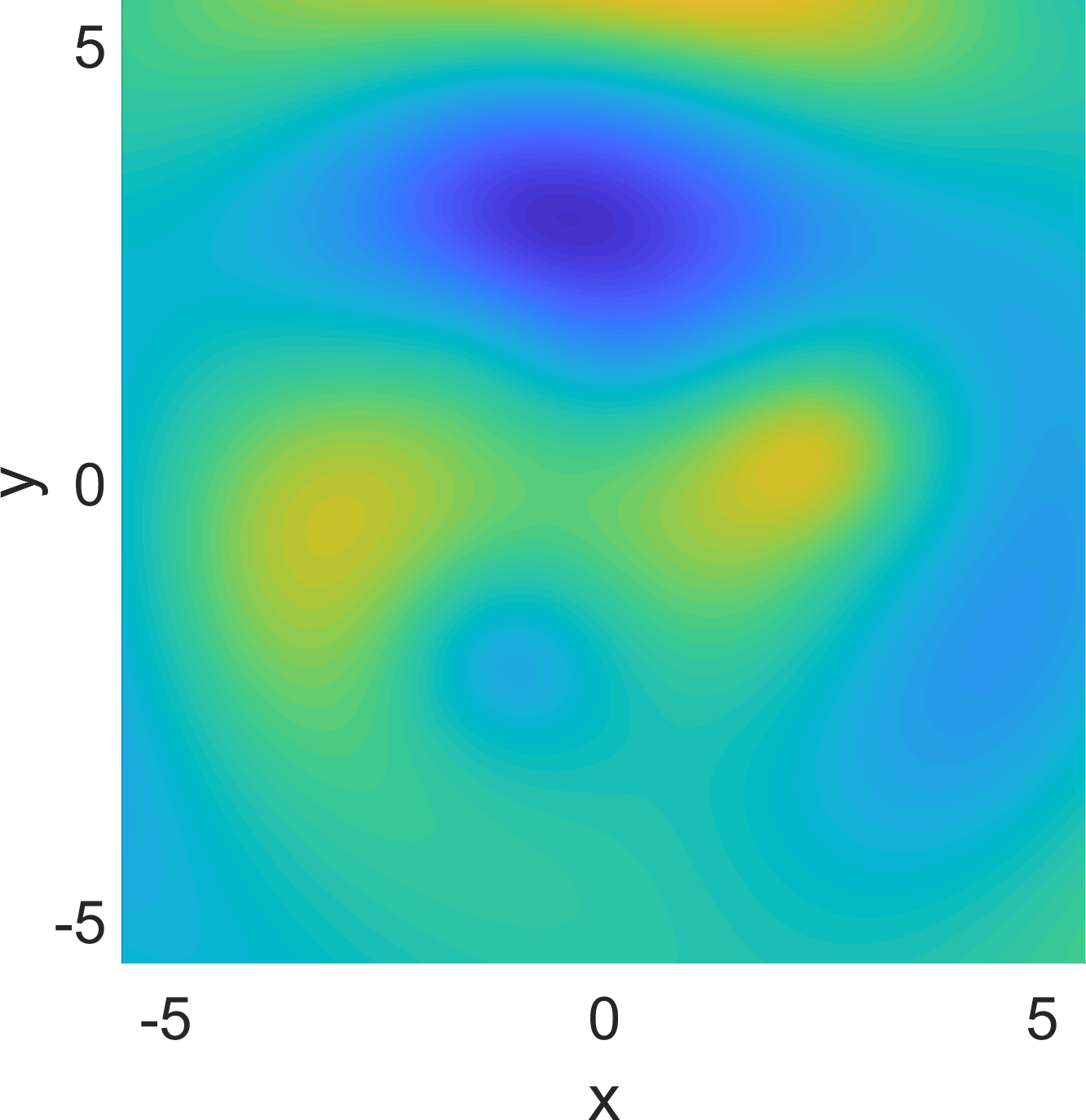}\hspace{2.5pt}
		\includegraphics[height = 0.33\linewidth]{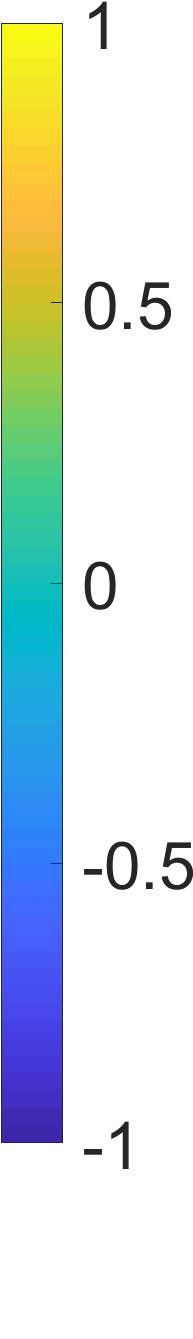}\\[2.5pt]
		\includegraphics[width = 0.49\linewidth]{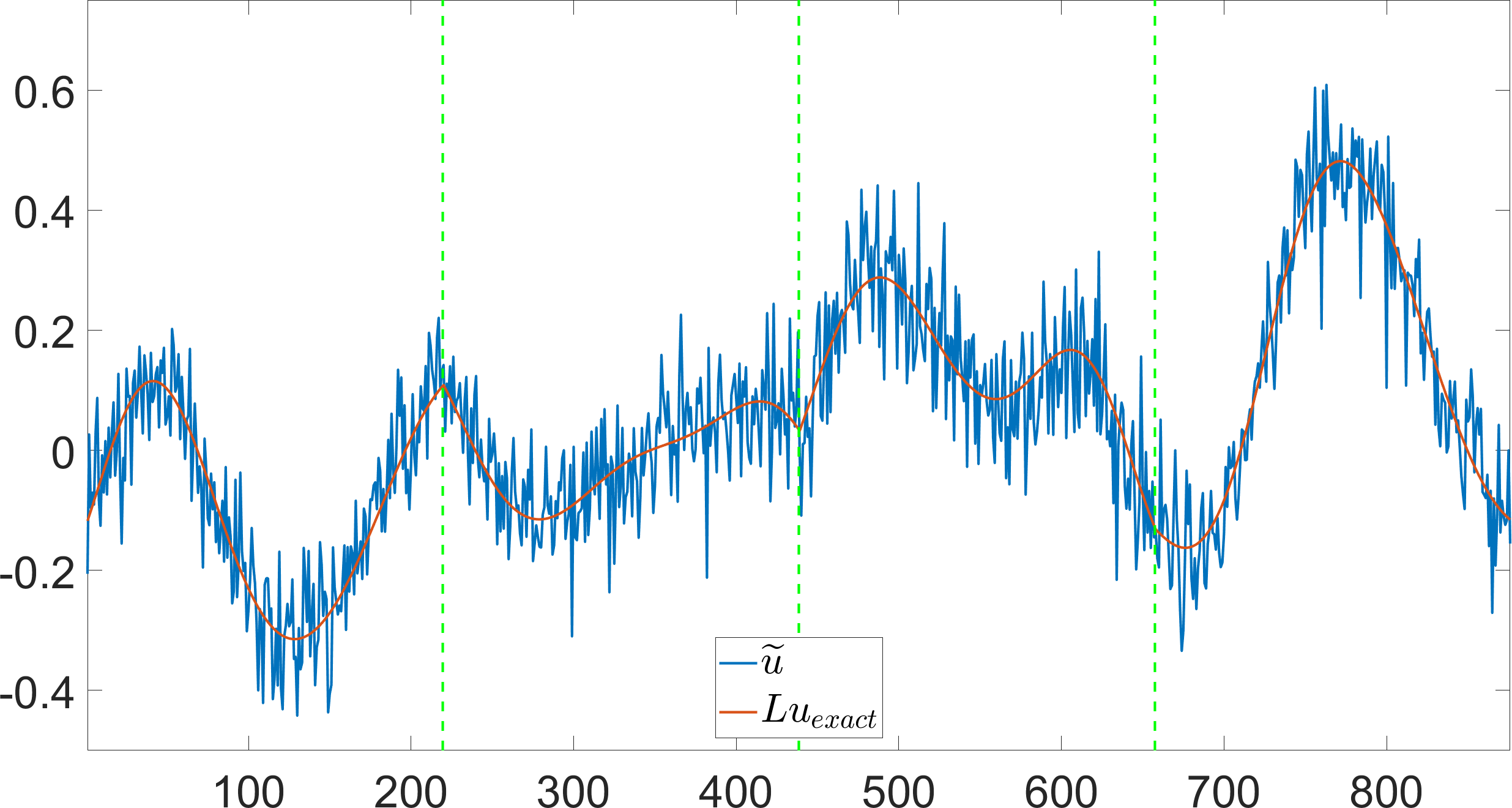}\hspace{2.5pt}
		\includegraphics[width = 0.49\linewidth]{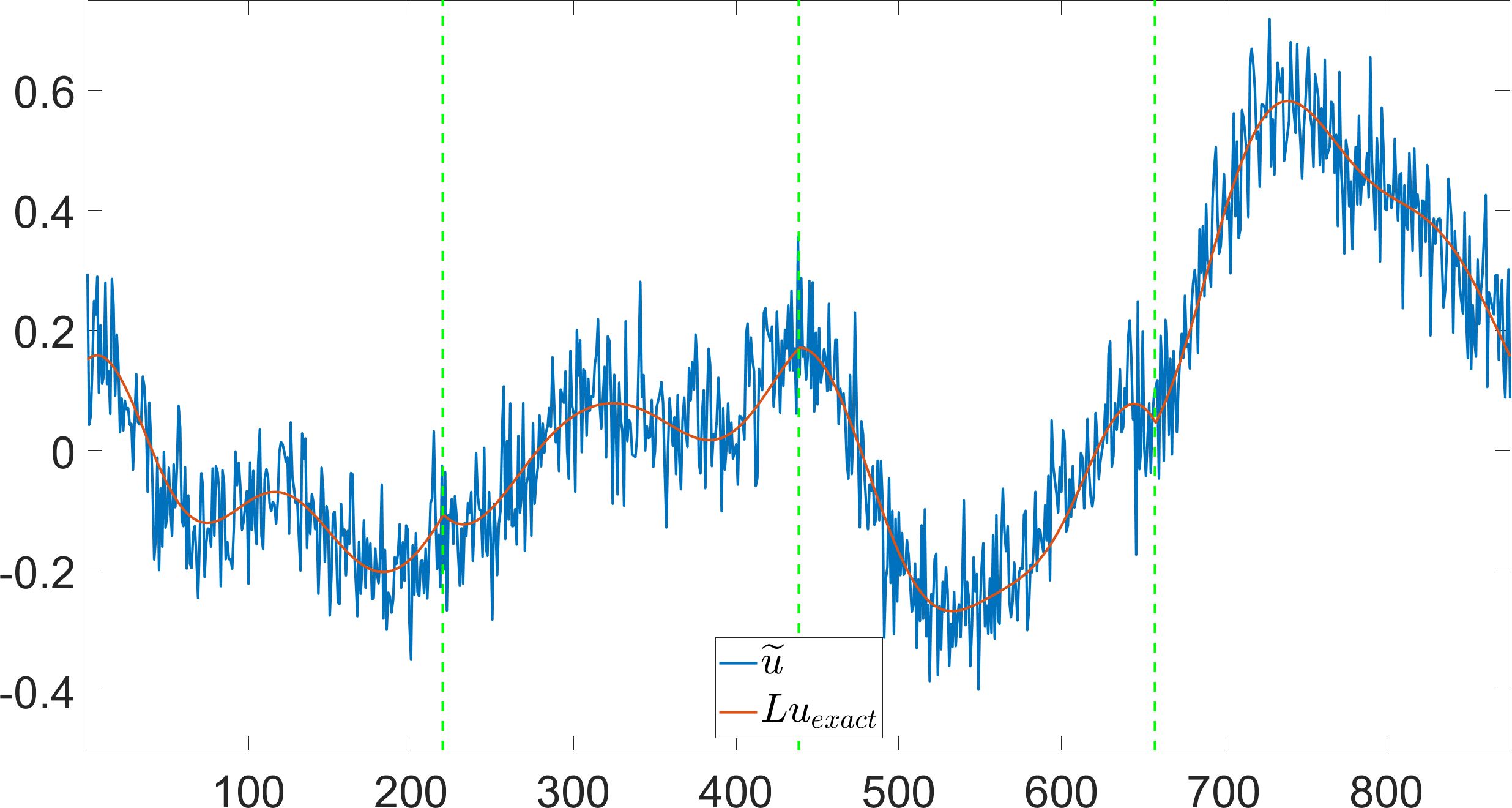}
		\caption{Top: real and complex part of the scattered wave for
				$\theta = 86.4^\circ$ (left and right) on $\Omega$, including the 10 buffer points to
				where the measurements are taken. Bottom: real and complex values of
				the same scattered wave (left and right) at the measurement points.
				The measurements are ordered counter clockwise starting from
				the top left corner of the domain and the dashed lines correspond to the
				corners.}
		\label{fig:scatwave}
\end{figure}


\subsection{Results}
\begin{figure}
		\centering
		\includegraphics[width = 0.49\linewidth]{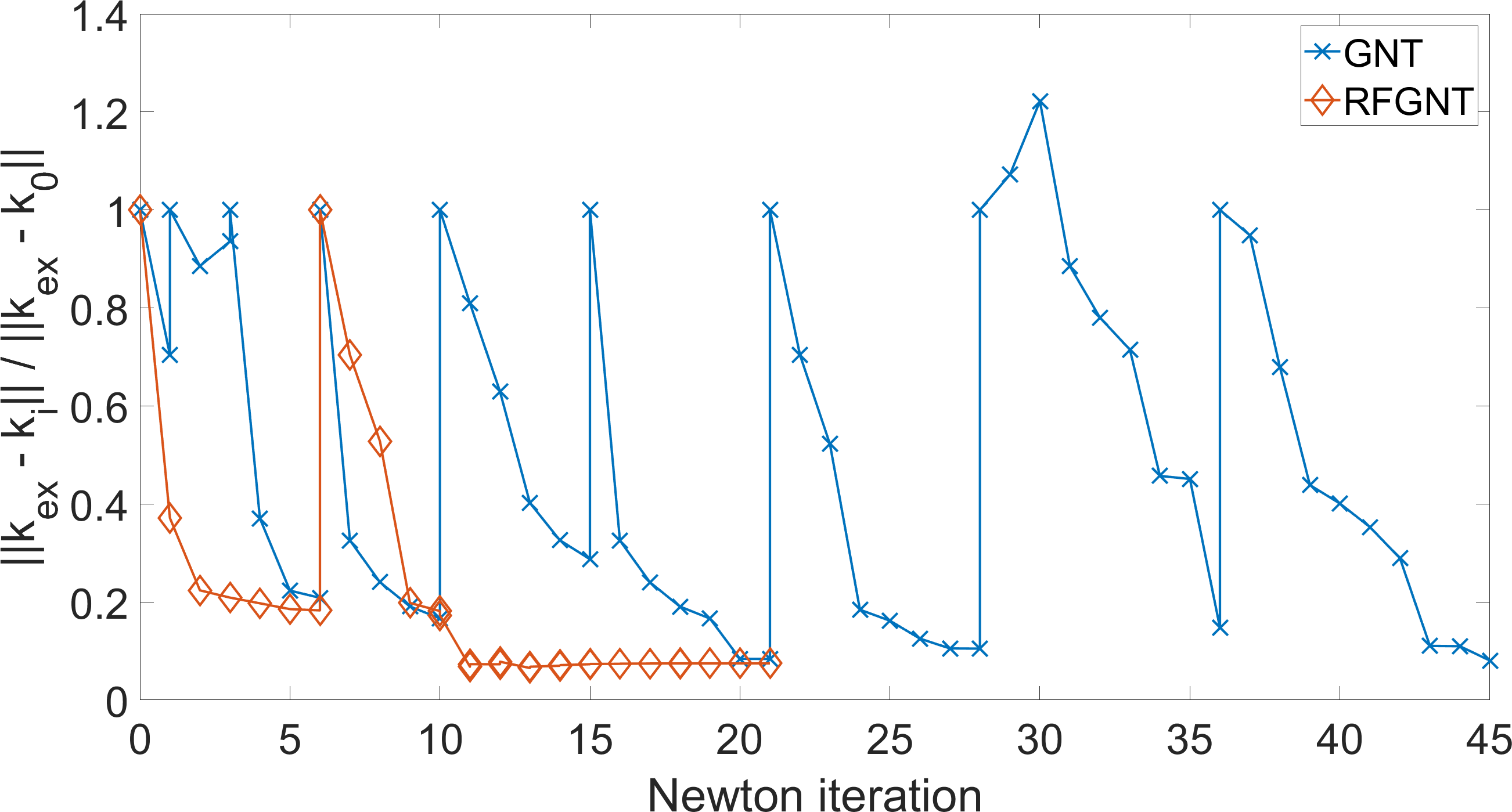}\hspace{2.5pt}
		\includegraphics[width = 0.49\linewidth]{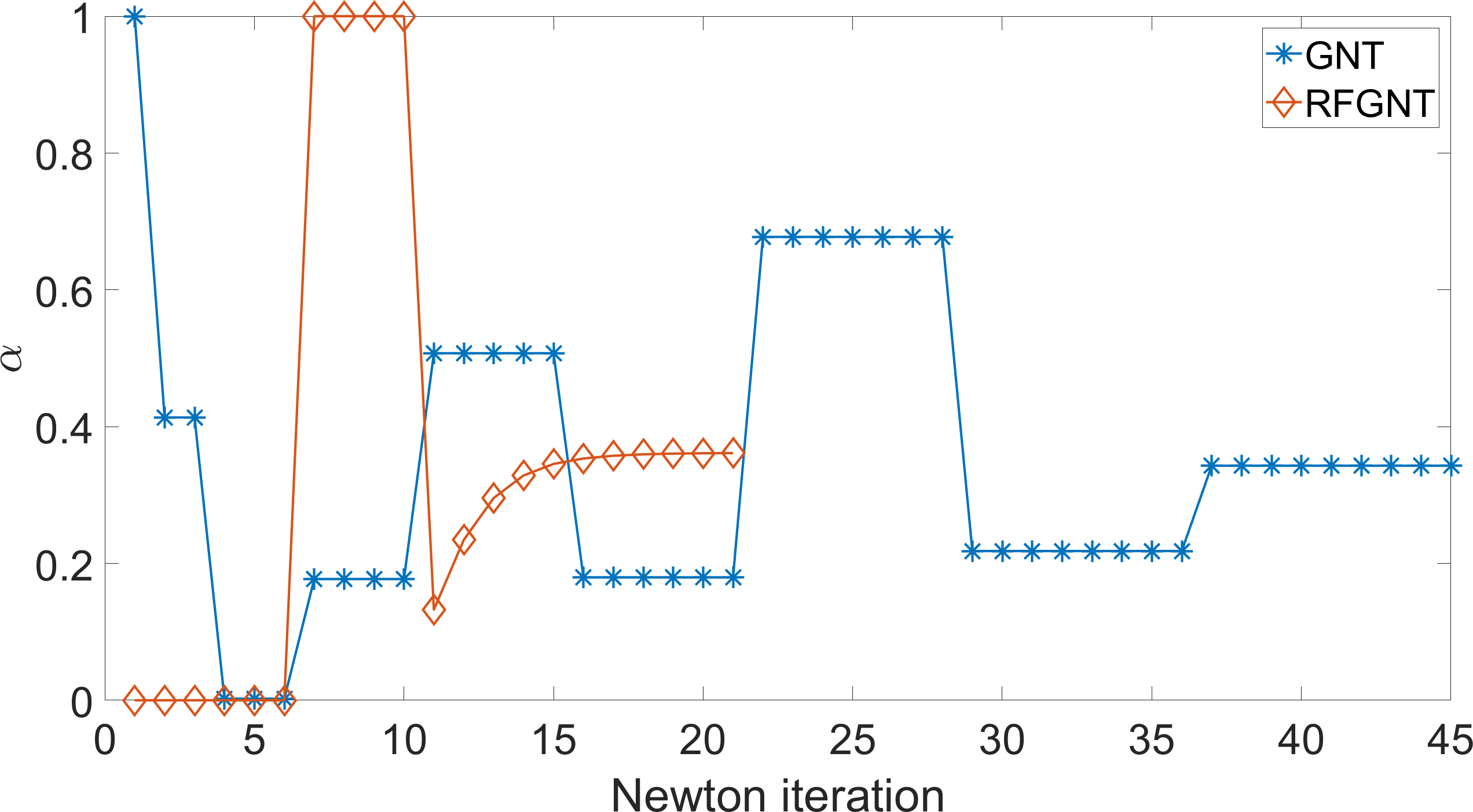}\\[2.5pt]
		\includegraphics[width = 0.49\linewidth]{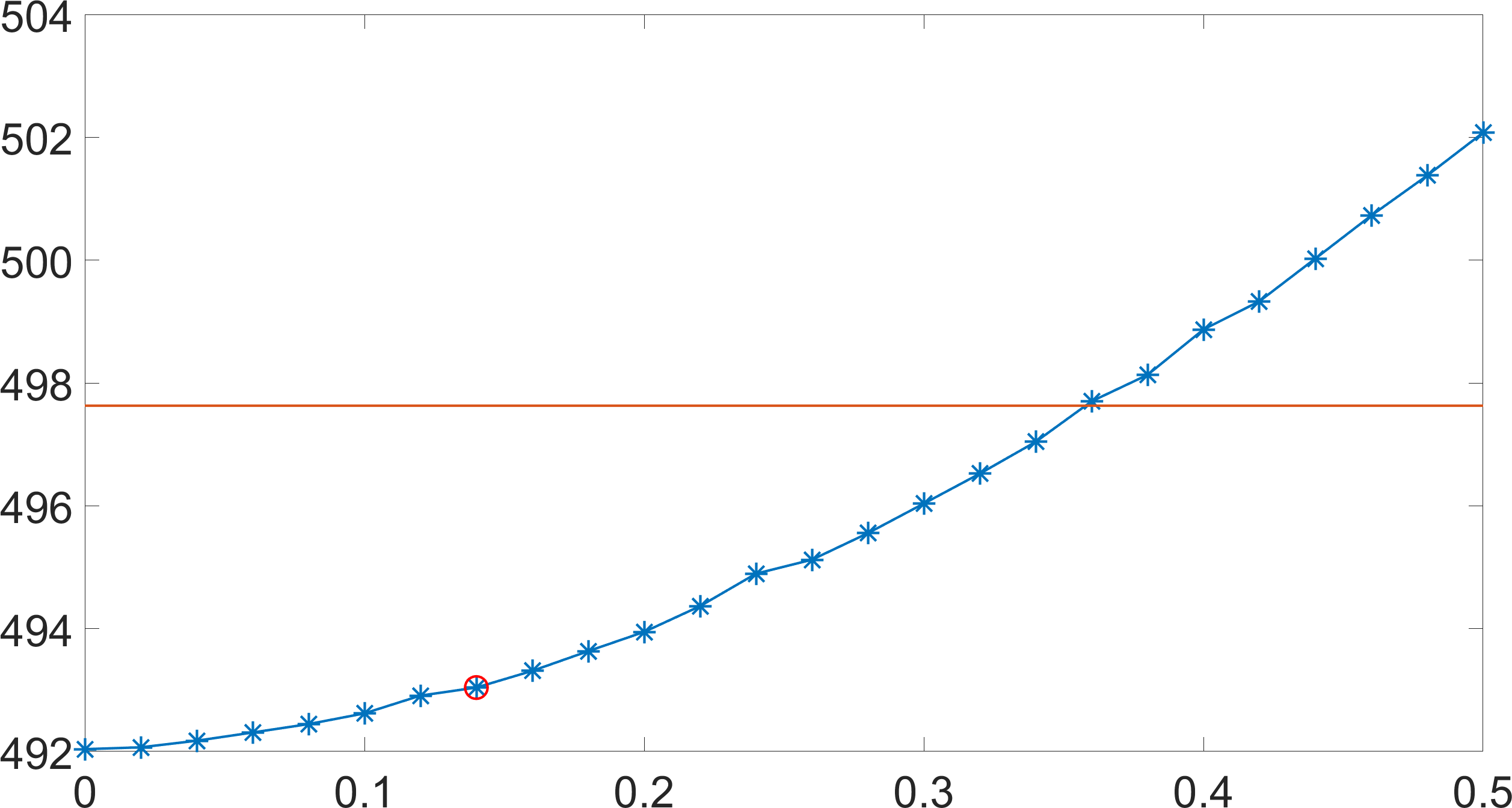}\hspace{2.5pt}
		\includegraphics[width = 0.49\linewidth]{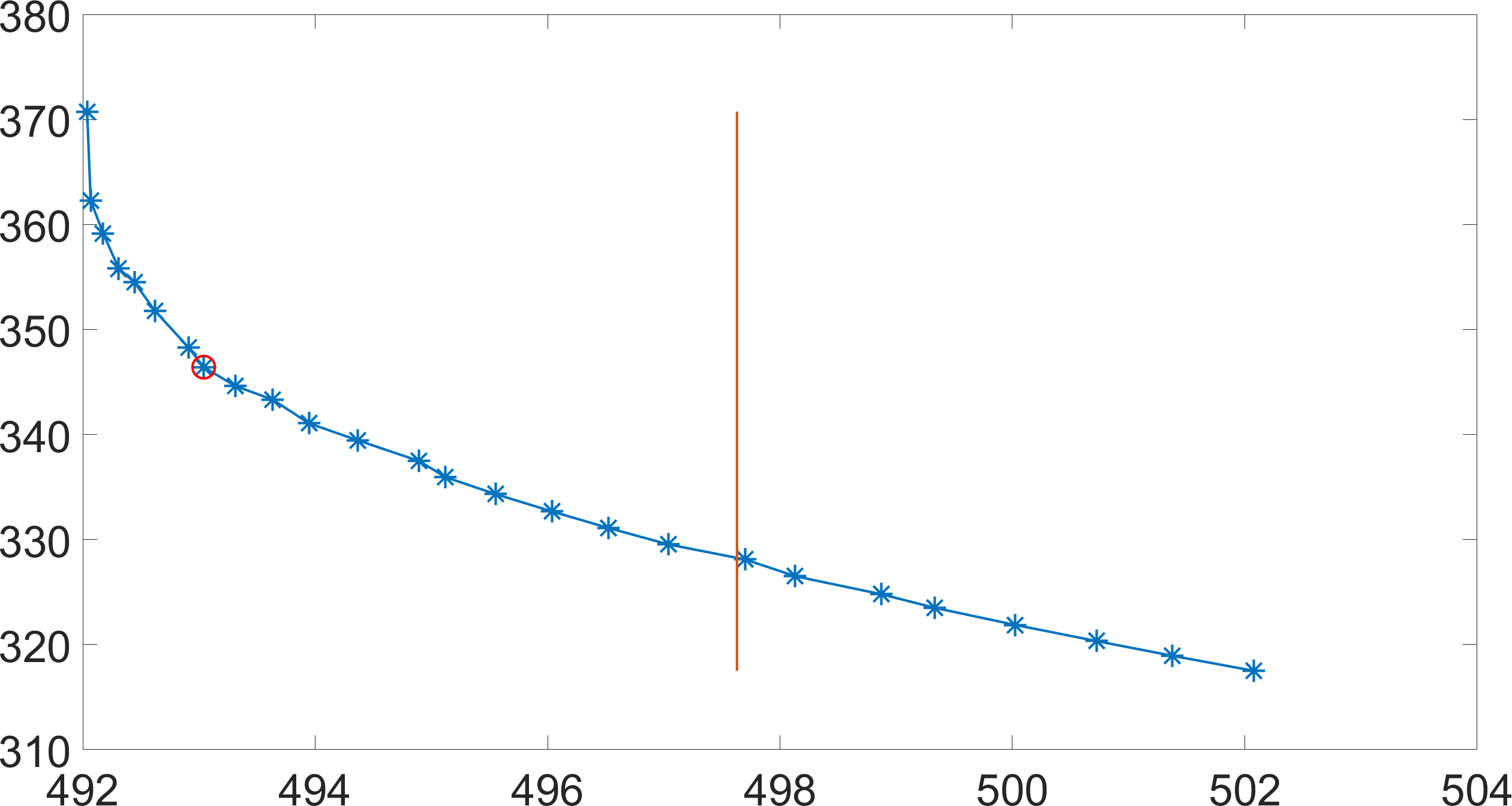}
		\caption{Top left: the relative error of GNT and RFGNT in each Newton iteration.
				We can clearly see at which point GNT restarts its Newton iterations.
				Top right: The regularization parameter used in each Newton iteration.
				Bottom: D-curve and L-curve. The red line corresponds to the discrepancy
				value $\eta\varepsilon$ and the point with the lowest error w.r.t. the exact
				solution is marked with a red circle.}
		\label{fig:details}
\end{figure}

In order to test the GNT and RFGNT algorithms, we consider two other reconstruction
methods. The first is Newton's method without a regularization term, but with the
discrepancy as an early stopping criterion. The second is the L-curve approach.
We calculated 26 points on the L-curve for $\alpha\in[0, 0.5]$ and selected the one
with the smallest error with respect to the exact wave number. Since the D-curve is
another way of looking at the L-curve, we can also use these points to see whether
or not the regularization parameter found by GNT and RFGNT is correct. Also, for GNT
and RFGNT we solved the problem using a backward (B), a central (C) and a forward (F)
finite difference approximation for the Hessian matrix vector product in the CG
iterations. However, since there was little difference in the quality of the
reconstructions, the early stopping solution and L-curve solution were calculated only
for the forward finite difference approximation. The details of all the reconstructions
are listed in \hypref{table}{tab:results} and \hypref{figures}{fig:recs}
and \ref{fig:details} illustrate some of these results for the forward finite
difference scheme.

When comparing the different reconstructions, we see that there is little difference
between the finite difference schemes, except in the number of gradient evaluations
(and hence the number of PDE, adjoint PDE and Helmholtz solves). We also see that while
GNT finds a good solution for the inverse problem, it requires a lot more Newton iterations
to do so. This can also be seen when looking at the relative error and regularization 
parameter in \hypref{figure}{fig:details}. Each time the Newton method is
restarted, progress is lost. While this approach was natural for linear problems
combined with a Krylov subspace method, it is inefficient for nonlinear problems.
The adaptations we made in order to derive the RFGNT method on the other hand seem to
be very effective. The method first needs to solve the problem for our initial choices
for $\alpha_0 = 0$ and $\alpha_1 = 1$, but due to the update of the initial estimate
for the Newton iterations, only 1 or 2 Newton iterations are needed afterwards for
\hypref{line}{alg:rfgnt:8} of \hypref{algorithm}{alg:rfgnt}. The convergence of the
regularization parameter is also drastically increased and by looking at the
discrepancy curve in \hypref{figure}{fig:details}, we see that the method
does indeed converge to the desired value for the regularization parameter.

When we compare GNT and RFGNT with early stopping, we see that the latter needs
less iterations before satisfying its stopping criterion. Then again, the early
stopping reconstruction has a much larger relative error than the other the
reconstructions. Using the L-curve approach, on the other hand, we find a
reconstruction that has the same quality as the GNT and RFGNT reconstructions.
However, while calculating a single point on the L-curve is cheaper than solving
the inverse problem with GNT or RFGNT, calculating all points the L-curve (26 in this case)
is much less efficient. This clearly illustrates the effectiveness off the automatic
regularization approach of RFGNT.

\begin{figure}
		\centering
		\includegraphics[width = 0.49\linewidth]{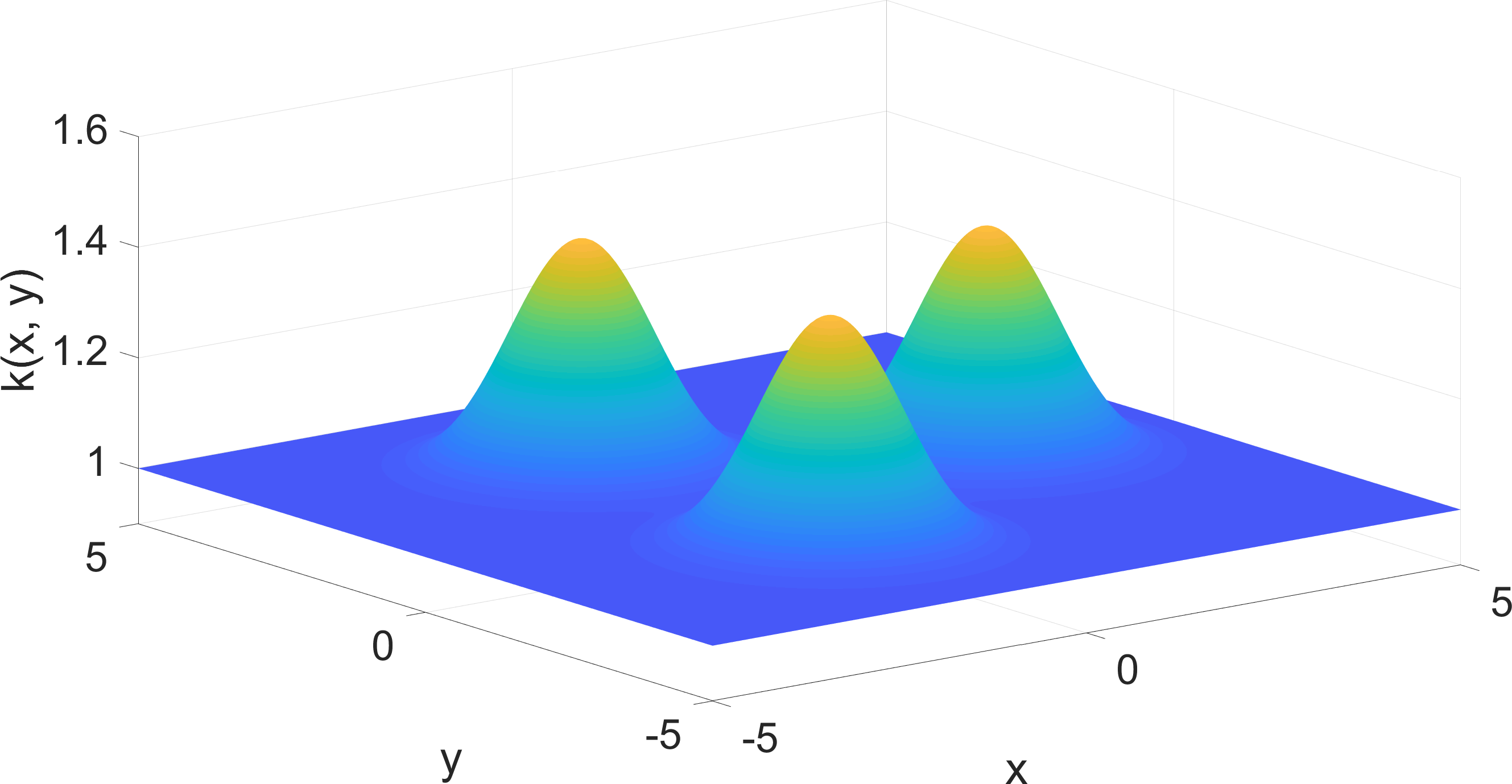}\\[2.5pt]
		\includegraphics[width = 0.49\linewidth]{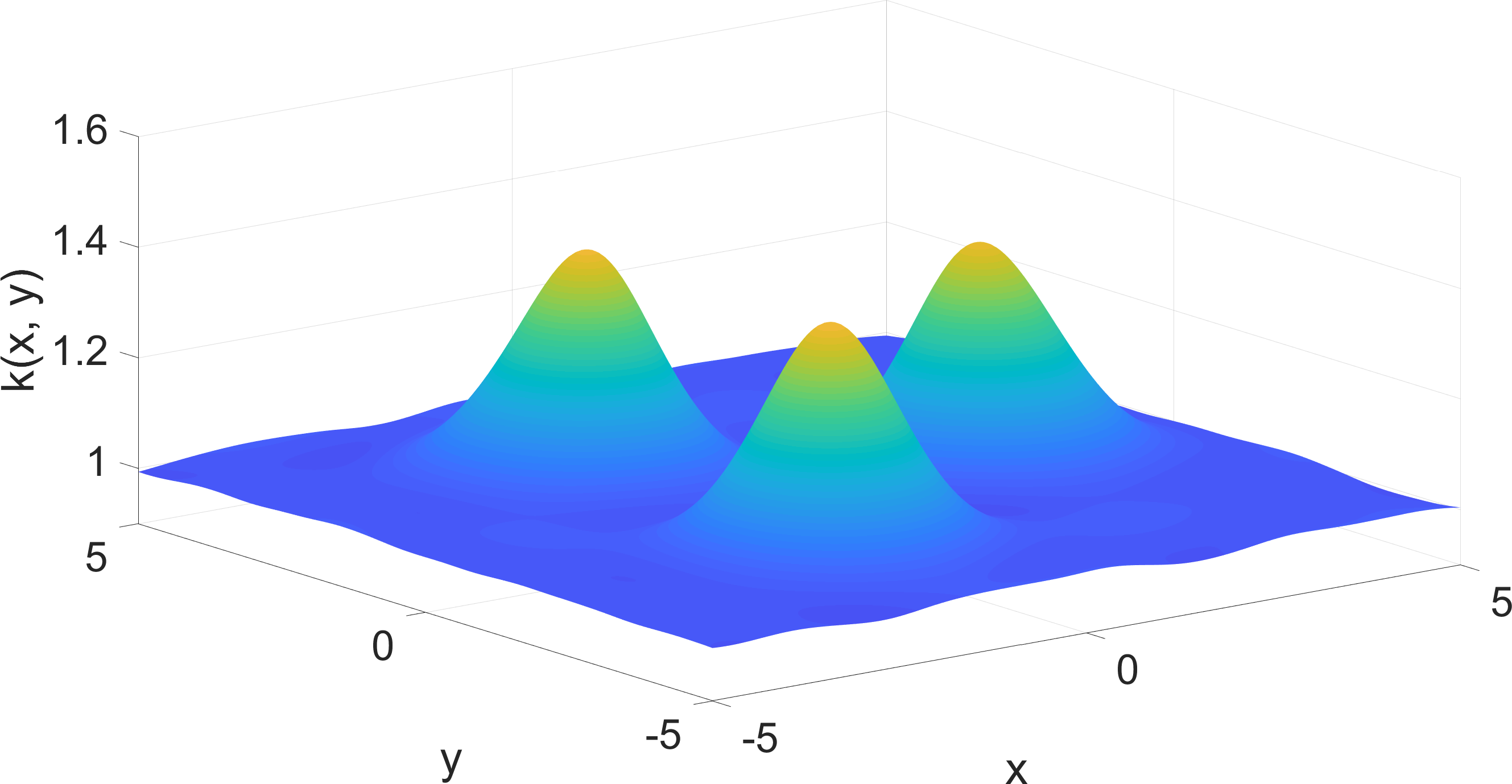}\hspace{2.5pt}
		\includegraphics[width = 0.49\linewidth]{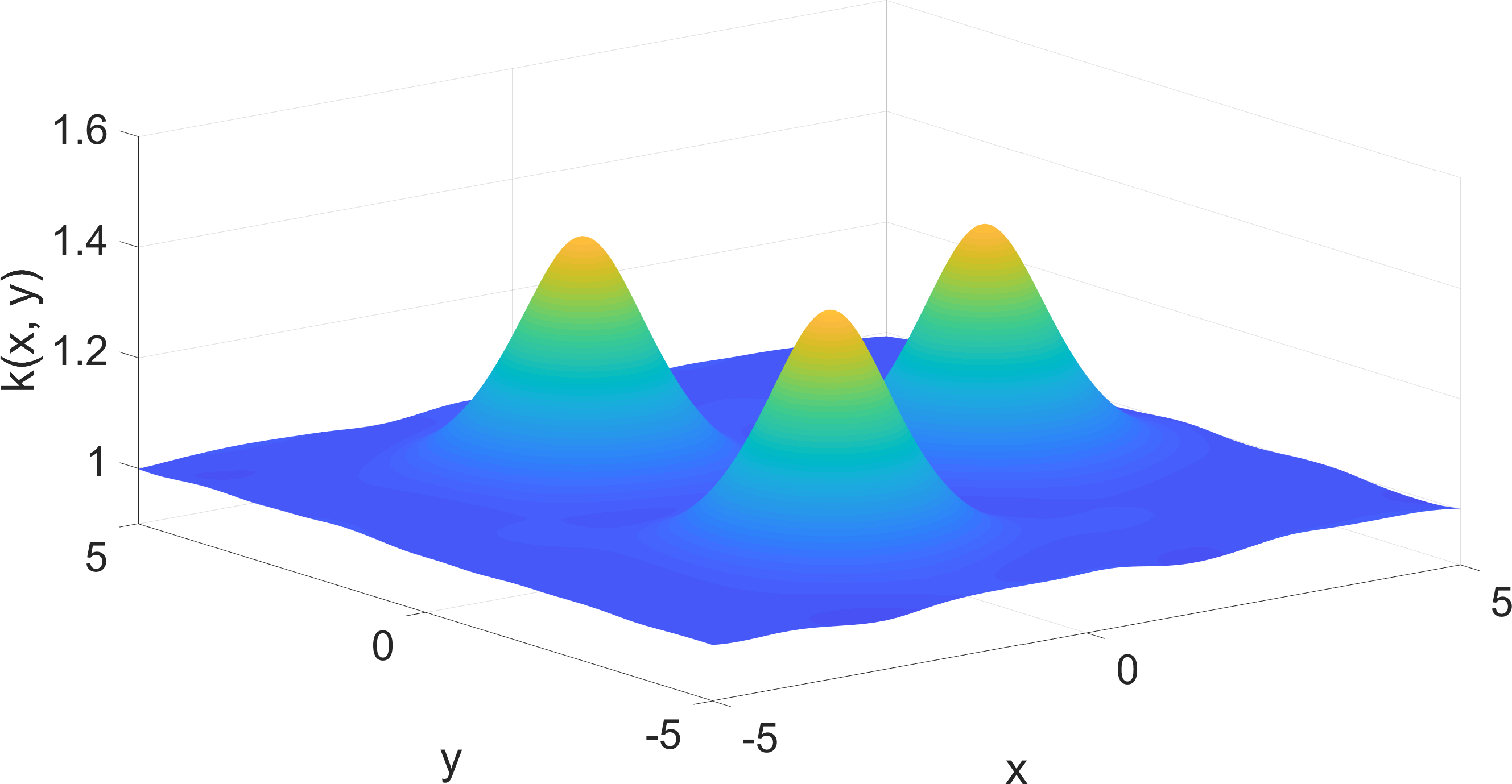}
		\includegraphics[width = 0.49\linewidth]{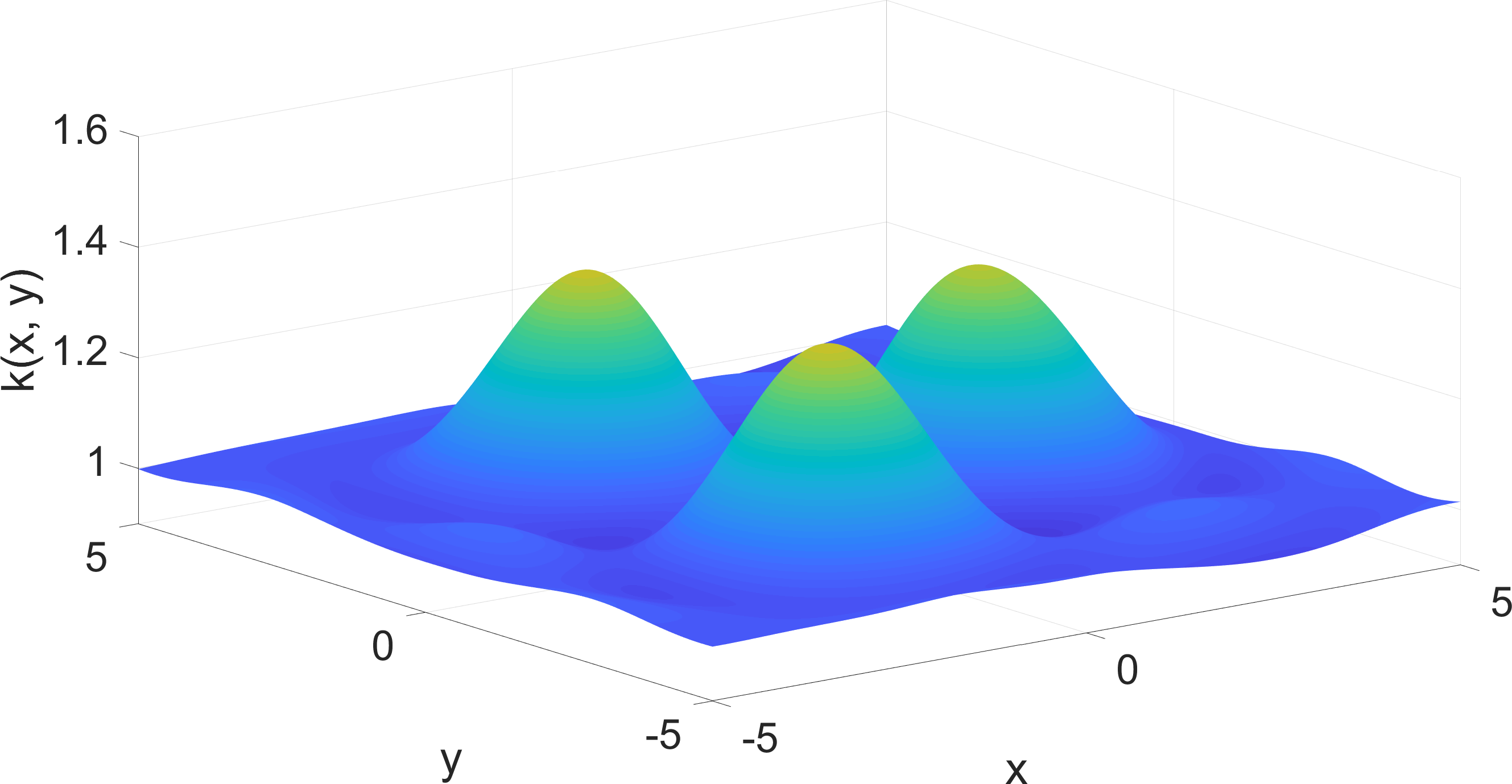}\hspace{2.5pt}
		\includegraphics[width = 0.49\linewidth]{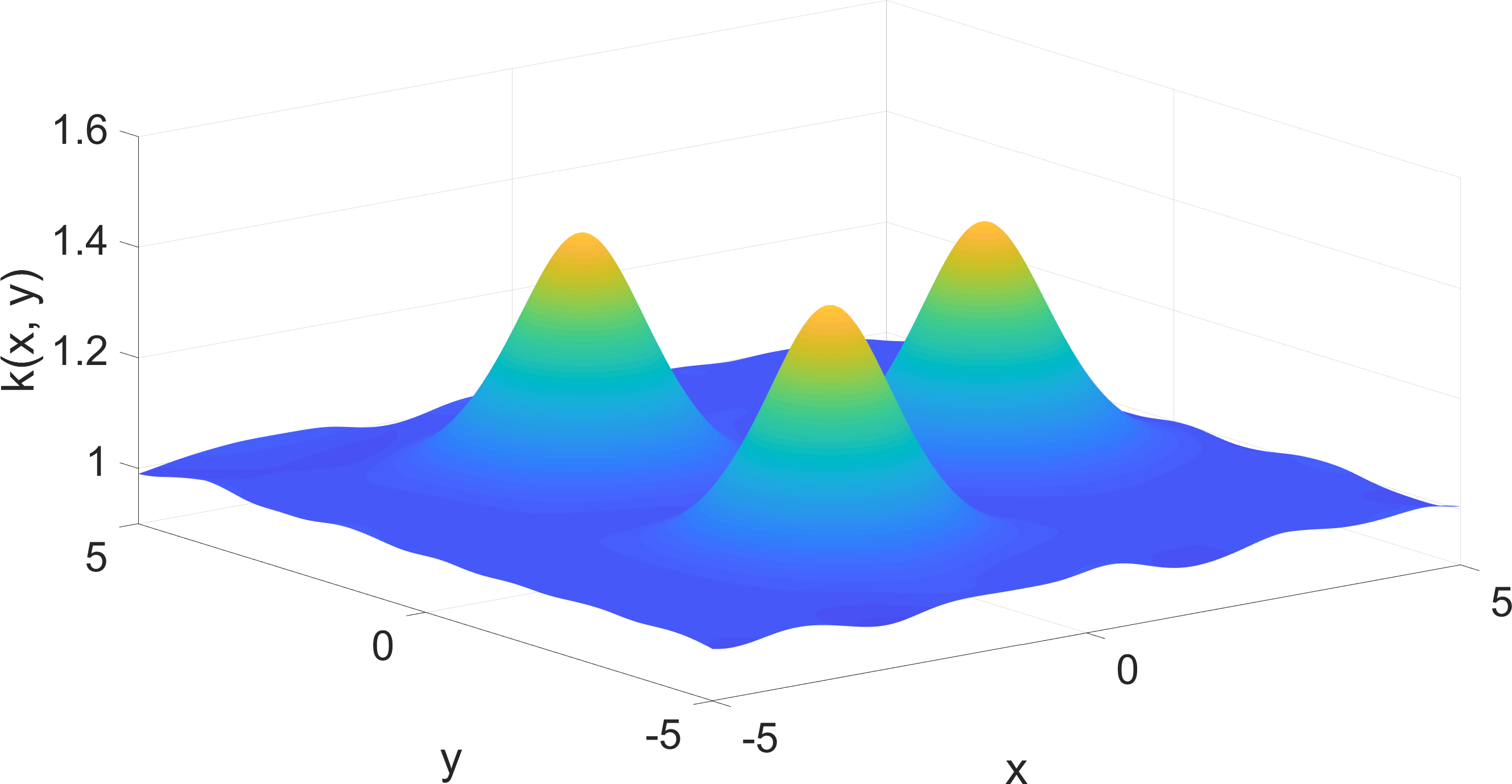}
		\caption{Top: the exact wave number $k$. Middle left: GNT reconstruction.
				Middle right: RFGNT reconstruction. Bottom left: early stopping reconstruction.
				Bottom right: L-curve reconstruction.}
		\label{fig:recs}
\end{figure}

\begin{table}
		\centering
		\scalebox{0.89}{
				\begin{tabular}{l|l||c|c|c|c|c|c|c|c|c}
						\multicolumn{2}{l||}{} & \rot{Newton iterations}& \rot{CG iterations} & \rot{$\mc{J}$ evaluations} & \rot{$\nabla\mc{J}$ evaluations}
								& \rot{PDE solves} & \rot{Adjoint solves} & \rot{$(\Delta + k^2)$  solves} & \rot{Relative error} & \rot{$\alpha$} \\ \hline\hline
						\multirow[c]{2}{*}{B} & GNT           & 54  & 175 & 101 & 108 & 209  & 108 & 15850 & 0.0803  & 0.3431  \\ \cline{2-11} 
																	& RFGNT 				& 21  & 73  & 44  & 42  & 86   & 42  & 6400  & 0.0746  & 0.3613  \\ \hline\hline
						\multirow[c]{2}{*}{C} & GNT           & 65  & 225 & 118 & 195 & 313  & 195 & 25400 & 0.0702  & 0.3398  \\ \cline{2-11} 
																	& RFGNT 				& 21  & 72  & 44  & 63  & 107  & 63  & 8500  & 0.0746  & 0.3613  \\ \hline\hline
						\multirow[c]{4}{*}{F} & GNT           & 54  & 175 & 101 & 108 & 209  & 108 & 15850 & 0.0804  & 0.3429  \\ \cline{2-11} 
																	& RFGNT         & 21  & 81  & 44  & 42  & 86   & 42  & 6400  & 0.0747  & 0.3613  \\ \cline{2-11}
																	& Early stopping & 3   & 7   & 8   & 6   & 14   & 6   & 1000  & 0.2089  & $\cdot$ \\ \cline{2-11}
																	& L-curve       & 9   & 37  & 20  & 18  & 38   & 18  & 2800  & 0.0871  & 0.1400  \\ \cline{2-11}
																	& L-curve (all) & 267 & 983 & 633 & 534 & 1167 & 534 & 85050 & $\cdot$ & $\cdot$
				\end{tabular}
		}
		\caption{Details from the different reconstructions for the different reconstruction
				methods and finite difference approximations for the Hessian matrix vector
				product in the CG iterations. The number of gradient evaluations is 2 or 3
				times the number of Newton iterations (depending on the finite difference
				scheme) which is equal to the number of solves of the adjoint PDE. The
				number of PDE solves is sum of the number of cost function evaluations
				and the number of gradient evaluations. The number of Helmholtz system solves
				is 50 times, i.e. the number of projection angles, the number of
				PDE and Adjoint solves.}
		\label{tab:results}
\end{table}


\section{Conclusion and remarks}
In this paper we describe two methods to solve a nonlinear inverse problem that
iteratively determine the solution and the regularization parameter. The first method,
generalized Newton-Tikhonov (GNT), is a direct generalization of the generalized
Arnoldi-Tikhonov method to nonlinear problems. However, this method turns out to
have a number of drawbacks that were not present in the original algorithm for linear
problems. In order to improve the method, we proposed the regula falsi generalized
Newton-Tikhonov method (RFGNT). We replace the secant update step from GNT with a
regula falsi approach and updating the initial guess for the Newton iterations with
every update of the regularization parameter. This decreases the number of Newton
iterations needed and finds a better value for the regularization parameter. Our
numerical experiments also show that this is computationally much more efficient than,
for example, calculating the L-curve or other grid based approaches to determine
the regularization parameter.

It should also be noted that in this paper we solve the PDE and the adjoint PDE
sequentially. By contrast, it is also possible to solve both simultaneously by
considering the Karush-Kuhn-Tucker conditions and using Newton's method combined
with a suitable preconditioner to find the optimum, see for example \cite{biros2005,
haber2000, haber2001, mardal2017}. The difficulty with this approach is finding a
suitable preconditioner for the problem, but it can easily be combined with the
proposed RFGNT method. This is because using this approach to solve the problem
for a fixed value of the regularization parameter $\alpha$ corresponds to
calculating $k_\alpha$ and $\mc{D}(\alpha)$, i.e. \hypref{lines}{alg:rfgnt:1},
\ref{alg:rfgnt:2} and \ref{alg:rfgnt:8} of \hypref{algorithm}{alg:rfgnt}.

Furthermore, although we used an inverse scattering problem as a test problem to
demonstrate the methods, no specific properties of this problem were used to derive
the methods and they are likely to be effective in many other inverse problems.
Future work therefore includes more in-depth analysis of the robustness of the methods
and its use in other application. We did, for example, not use any preconditioner
for the solution of the Helmholtz problem or the inner CG iterations. However,
there has been done many interesting work concerning preconditioners for the
Helmholtz equations using multigrid methods \cite{cools2013, erlangga2006}. Including
these in the optimization procedure could reduce the number of required solves
of the Helmholtz equation and CG iterations. Another possible issue is the fact
that since Newton can only be used for local minimization a proper initial
estimate needs to be determined.


\section*{Acknowledgement}
The authors wish to thank the Department of Mathematics and Computer Science,
University of Antwerp, for financial support.


\section*{References}
\bibliographystyle{elsarticle-num}
\bibliography{References}

\begin{thebibliography}{10}
\expandafter\ifx\csname url\endcsname\relax
  \def\url#1{\texttt{#1}}\fi
\expandafter\ifx\csname urlprefix\endcsname\relax\def\urlprefix{URL }\fi
\expandafter\ifx\csname href\endcsname\relax
  \def\href#1#2{#2} \def\path#1{#1}\fi

\bibitem{abdoulaev2005}
G.~S. Abdoulaev, K.~Ren, A.~H. Hielscher, Optical tomography as a
  {PDE}-constrained optimization problem, Inverse Problems 21~(5) (2005) 1507.

\bibitem{bruckner2017}
F.~Bruckner, C.~Abert, G.~Wautischer, C.~Huber, C.~Vogler, M.~Hinze, D.~Suess,
  Solving large-scale inverse magnetostatic problems using the adjoint method,
  Scientific reports 7 (2017) 40816.

\bibitem{kaebe2009}
C.~Kaebe, J.~H. Maruhn, E.~W. Sachs, Adjoint-based {M}onte {C}arlo calibration
  of financial market models, Finance and Stochastics 13~(3) (2009) 351--379.

\bibitem{inthout2010}
K.~J. {In't Hout}, S.~Foulon, {ADI} finite difference schemes for option
  pricing in the {H}eston model with correlation, International journal of
  numerical analysis and modeling 7~(2) (2010) 303--320.

\bibitem{nocedal2006}
J.~Nocedal, S.~J. Wright, Numerical optimization, Springer, 2006.

\bibitem{tortorelli1994}
D.~A. Tortorelli, P.~Michaleris, Design sensitivity analysis: overview and
  review, Inverse problems in Engineering 1~(1) (1994) 71--105.

\bibitem{plessix2006}
R.~Plessix, A review of the adjoint-state method for computing the gradient of
  a functional with geophysical applications, Geophysical Journal International
  167~(2) (2006) 495--503.

\bibitem{jadamba2017}
B.~Jadamba, A.~A. Khan, A.~A. Oberai, M.~Sama, First-order and second-order
  adjoint methods for parameter identification problems with an application to
  the elasticity imaging inverse problem, Inverse Problems in Science and
  Engineering 25~(12) (2017) 1768--1787.

\bibitem{calvetti1999}
D.~Calvetti, G.~H. Golub, L.~Reichel, Estimation of the {L}-curve via {L}anczos
  bidiagonalization, BIT Numerical Mathematics 39~(4) (1999) 603--619.

\bibitem{calvetti2004}
D.~Calvetti, L.~Reichel, A.~Shuibi, L-curve and curvature bounds for tikhonov
  regularization, Numerical Algorithms 35~(2--4) (2004) 301--314.

\bibitem{hansen2010}
P.~C. Hansen, Discrete inverse problems: insight and algorithms, Vol.~7, Siam,
  2010.

\bibitem{vogel2002}
C.~R. Vogel, Computational methods for inverse problems, SIAM, 2002.

\bibitem{gazzola2014_2}
S.~Gazzola, J.~G. Nagy, Generalized {A}rnoldi-{T}ikhonov method for sparse
  reconstruction, SIAM Journal on Scientific Computing 36~(2) (2014)
  B225--B247.

\bibitem{gazzola2014}
S.~Gazzola, P.~Novati, Automatic parameter setting for {A}rnoldi-{T}ikhonov
  methods, Journal of Computational and Applied Mathematics 256 (2014)
  180--195.

\bibitem{gazzola2015}
S.~Gazzola, P.~Novati, M.~R. Russo, On krylov projection methods and tikhonov
  regularization, Electron. Trans. Numer. Anal 44 (2015) 83--123.

\bibitem{hansen1992}
P.~C. Hansen, Analysis of discrete ill-posed problems by means of the
  {L}-curve, SIAM review 34~(4) (1992) 561--580.

\bibitem{morozov1984}
V.~A. Morozov, Methods for Solving Incorrectly Posed Problems,
  Springer-Verslag, 1984.

\bibitem{saad2003}
Y.~Saad, Iterative methods for sparse linear systems, Vol.~82, siam, 2003.

\bibitem{vandervorst2003}
H.~A. {Van der Vorst}, Iterative Krylov methods for large linear systems,
  Vol.~13, Cambridge University Press, 2003.

\bibitem{saad1986}
Y.~Saad, M.~H. Schultz, {GMRES}: A generalized minimal residual algorithm for
  solving nonsymmetric linear systems, SIAM Journal on scientific and
  statistical computing 7~(3) (1986) 856--869.

\bibitem{shewchuk1994}
J.~R. Shewchuk, et~al., An introduction to the conjugate gradient method
  without the agonizing pain (1994).

\bibitem{wang1992}
Z.~Wang, I.~M. Navon, F.~X.~L. Dimet, X.~Zou, The second order adjoint
  analysis: theory and applications, Meteorology and atmospheric physics
  50~(1-3) (1992) 3--20.

\bibitem{mccurdy2004}
C.~W. McCurdy, M.~Baertschy, T.~N. Rescigno, Solving the three-body coulomb
  breakup problem using exterior complex scaling, Journal of Physics B: Atomic,
  Molecular and Optical Physics 37~(17) (2004) R137.

\bibitem{biros2005}
G.~Biros, O.~Ghattas, Parallel {L}agrange--{N}ewton--{K}rylov--{S}chur methods
  for {PDE}-constrained optimization. part {I}: The {K}rylov--{S}chur solver,
  SIAM Journal on Scientific Computing 27~(2) (2005) 687--713.

\bibitem{haber2000}
E.~Haber, U.~M. Ascher, D.~Oldenburg, On optimization techniques for solving
  nonlinear inverse problems, Inverse Problems 16~(5) (2000) 1263.

\bibitem{haber2001}
E.~Haber, U.~M. Ascher, Preconditioned all-at-once methods for large, sparse
  parameter estimation problems, Inverse Problems 17~(6) (2001) 1847.

\bibitem{mardal2017}
K.-A. Mardal, B.~F. Nielsen, M.~Nordaas, Robust preconditioners for
  {PDE}-constrained optimization with limited observations, BIT Numerical
  Mathematics 57~(2) (2017) 405--431.

\bibitem{cools2013}
S.~Cools, W.~Vanroose, Local {F}ourier analysis of the complex shifted
  {L}aplacian preconditioner for {H}elmholtz problems, Numerical Linear Algebra
  with Applications 20~(4) (2013) 575--597.

\bibitem{erlangga2006}
Y.~A. Erlangga, C.~Vuik, C.~W. Oosterlee, Comparison of multigrid and
  incomplete {LU} shifted-{L}aplace preconditioners for the inhomogeneous
  {H}elmholtz equation, Applied numerical mathematics 56~(5) (2006) 648--666.

\end{thebibliography}


\end{document}